\documentclass[11pt]{article}

\usepackage{amsmath}
\usepackage{amsfonts}
\usepackage{amssymb}

\makeatletter

\@addtoreset{equation}{section}
\makeatother 

\begin{document}

\newfont{\cyr}{wncyr10 scaled\magstep1}
\newcommand{\sh}{{\cyr\char '130}}
\newcommand{\sha}{\mbox{\sh}}

\newtheorem{Theorem}{Theorem}[subsection]
\newtheorem{Conjecture}[Theorem]{Conjecture}
\newtheorem{Proposition}[Theorem]{Proposition}
\newtheorem{Lemma}[Theorem]{Lemma}
\newtheorem{Corollary}[Theorem]{Corollary}
\newtheorem{Remark}[Theorem]{Remark}
\newtheorem{Example}[Theorem]{Example}
\newtheorem{Definition}[Theorem]{Definition}
\newtheorem{Question}[Theorem]{Question}
\newtheorem{TTheorem}{Theorem} 
\def\theTTheorem{\Alph{TTheorem}}

\newcommand{\ord}{\mathop{\rm ord}\nolimits}
\newcommand{\ddiv}{\mathop{\rm div}\nolimits}

\newcommand{\mm}{\mathop{\rm mod}\nolimits}
\newcommand{\Spec}{\mathop{\rm Spec}\nolimits}
\newcommand{\gr}{\mathop{\rm gr}\nolimits}
\newcommand{\tr}{\mathop{\rm tr}\nolimits} 
\newcommand{\Cor}{\mathop{\rm Cor}\nolimits}
\newcommand{\Coker}{\mathop{\rm Coker}\nolimits}
\newcommand{\Image}{\mathop{\rm Image}\nolimits}
\newcommand{\Aut}{\mathop{\rm Aut}\nolimits}
\newcommand{\tors}{\mathop{\rm tors}\nolimits}
\newcommand{\Res}{\mathop{\rm Res}\nolimits}
\newcommand{\Gal}{\mathop{\rm Gal}\nolimits}
\newcommand{\Div}{\mathop{\rm Div}\nolimits}
\newcommand{\Tor}{\mathop{\rm Tor}\nolimits}
\newcommand{\Fitt}{\mathop{\rm Fitt}\nolimits}
\newcommand{\Ann}{\mathop{\rm Ann}\nolimits}
\newcommand{\ii}{\mathop{\rm Image}\nolimits}
\newcommand{\re}{\mathop{\rm Re}\nolimits}
\newcommand{\Sel}{\mathop{\rm Sel}\nolimits}
\newcommand{\length}{\mathop{\rm length}\nolimits}
\newcommand{\cc}{\mathop{\rm char}\nolimits}
\newcommand{\Cl}{\mathop{\rm Cl}\nolimits}
\newcommand{\Gel}{\mathop{\mathfrak l}\nolimits}
\newcommand{\Gn}{\mathop{\mathfrak n}\nolimits}
\newcommand{\Gr}{\mathop{\mathfrak r}\nolimits}
\newcommand{\Gm}{\mathop{\mathfrak m}\nolimits}
\newcommand{\Gd}{\mathop{\mathfrak d}\nolimits}
\newcommand{\pr}{\mathop{\rm pr}\nolimits}
\newcommand{\F}{\mathop{\rm Frob}\nolimits}
\newcommand{\cd}{\mathop{\rm cd}\nolimits}
\newcommand{\Tam}{\mathop{\rm Tam}\nolimits}

\newcommand{\specialarrow}[3]{\mathrel{\mathop{#1}\limits^{#2}_{#3}}}
\newcommand{\superarrow}[2]{\mathrel{\mathop{#1}\limits_{#2}}}
\newcommand{\supperarrow}[2]{\mathrel{\mathop{#1}\limits^{#2}}}

\newcommand{\uetuki}[2]{{\mathop{#1}\limits^{#2}}}

\newcommand{\qed}
{{\unskip\nobreak\hfil\penalty50\quad\null\nobreak\hfil{Q.E.D.}\parfillskip0pt\finalhyphendemerits0\par\medskip}}

\newcommand{\pf}{\mbox{\it pf}}
\newcommand{\Pf}{\mbox{\it Pf}}
\newcommand{\llra}{\mathrel{- \hspace{-2mm} - \hspace{-2mm} \longrightarrow}}
\newcommand{\hhh}{\hspace*{3.5mm}}
\newcommand{\gl}{\mathop{\rm GL}\nolimits}
\newcommand{\www}{\textstyle\bigwedge\nolimits}
\newcommand{\s}{\mathop{\rm S}\nolimits}
\newcommand{\rank}{\mathop{\rm rank}\nolimits}
\newcommand{\image}{\mathop{\rm Im}\nolimits}
\def\rmn#1{\uppercase\expandafter{\romannumeral#1}}
\newcommand{\sgn}{\mathop{\rm sgn}\nolimits}
\newcommand{\Hom}{\mathop{\mbox{\rm Hom}}\nolimits}
\newcommand{\Ker}{\mathop{\mbox{\rm Ker}}\nolimits}
\newcommand{\Ext}{\mathop{\mbox{\rm Ext}}\nolimits}
\newcommand{\maru}{\mathrel{
\setlength{\unitlength}{.1cm}
\begin{picture}(1.5,2.5)(0,0)
\put(0.75,1){\circle{0.75}}
\end{picture}
}}

\newcommand{\mapdown}[1]{\Big\downarrow
 \rlap{$\vcenter{\hbox{$\scriptstyle#1$}}$ }}

\newcommand{\mapup}[1]{\Big\uparrow
 \rlap{$\vcenter{\hbox{$\scriptstyle#1$}}$ }}

\title{The structure of Selmer groups 
of elliptic curves and modular symbols}
\author{Masato {\sc Kurihara}}
\date{}
\maketitle

For an elliptic curve over the rational number field and 
a prime number $p$,   
we study the structure of the classical Selmer group of 
$p$-power torsion points. 
In our previous paper \cite{Ku6}, 
assuming the main conjecture and the non-degeneracy of 
the $p$-adic height pairing,
we proved that 
the structure of the Selmer group with respect to 
$p$-power torsion points
is determined by some analytic 
elements $\tilde{\delta}_{m}$ 
defined from 
modular symbols (see Theorem \ref{MT2} below). 
In this paper, we do not assume the main conjecture 
nor the non-degeneracy of the $p$-adic height pairing, and 
study the structure of Selmer groups 
(see Theorems \ref{MMT1} and \ref{MMT2}), 
using these analytic elements and Kolyvagin systems of 
Gauss sum type.

\addtolength{\textheight}{36pt}

\section{Introduction} \label{section1}

\subsection{Structure theorem of Selmer groups} \label{subsec11}

Let $E$ be an elliptic curve over ${\bf Q}$. 
Iwasawa theory, especially the main conjecture gives 
a formula on 
the order of the Tate Shafarevich group 
by using the $p$-adic $L$-function (cf. Schneider \cite{Sch0}). 
In this paper, as a sequel of \cite{Ku3}, \cite{Ku5} and \cite{Ku6}, 
we show that we can derive  
more information than the order,  
on the structure of the Selmer group and the Tate Shafarevich group 
from analytic quantities, in the setting of our paper, from modular symbols.

In this paper, we consider a prime number $p$ such that \\
(i) $p$ is a good ordinary prime $>2$ for $E$, \\
(ii) the action of $G_{{\bf Q}}$ on the Tate module 
$T_{p}(E)$ is surjective 
where $G_{{\bf Q}}$ is the absolute Galois group of ${\bf Q}$, \\
(iii) the (algebraic) $\mu$-invariant of $(E, {\bf Q}_{\infty}/{\bf Q})$ 
is zero where ${\bf Q}_{\infty}/{\bf Q}$ is the 
cyclotomic ${\bf Z}_{p}$-extension, 
namely the Selmer group $\Sel(E/{\bf Q}_{\infty}, E[p^{\infty}])$ 
(for the definition, see below) 
is a cofinitely generated ${\bf Z}_{p}$-module, \\
(iv) $p$ does not divide the Tamagawa factor 
$\Tam(E)=\Pi_{\ell:{\rm bad}} (E({\bf Q}_{\ell}):E^{0}({\bf Q}_{\ell}))$, 
and $p$ does not divide $\#E({\bf F}_{p})$ 
(namely not anomalous). \\

We note that the property (iii) is a conjecture of Greenberg 
since we are assuming (ii).

For a positive integer $N>0$, we denote by $E[p^{N}]$ the 
Galois module of $p^{N}$-torsion points, and 
$E[p^{\infty}]=\bigcup_{N>0} E[p^{N}]$.  
For an algebraic extension $F/{\bf Q}$, 
$\Sel(E/F, E[p^{N}])$ is the classical Selmer group 
defined by 
$$\Sel(E/F, E[p^{N}])=
\Ker(H^{1}(F, E[p^{N}]) \longrightarrow 
\prod_{v } H^{1}(F_{v}, E[p^{N}])
/E(F_{v}) \otimes {\bf Z}/p^{N}),$$
so $\Sel(E/F, E[p^{N}])$ sits in an exact sequence 
$$0 \longrightarrow E(F) \otimes {\bf Z}/p^{N} \longrightarrow 
\Sel(E/F, E[p^{N}]) \longrightarrow \sha(E/F)[p^{N}]
\longrightarrow 0$$
where $\sha(E/F)$ is the Tate Shafarevich group over $F$.  
We define $\Sel(E/F, E[p^{\infty}])= 
{\lim\limits_{\longrightarrow}}\Sel(E/F, E[p^{N}])$. 

Let ${\cal P}^{(N)}$ be the set of prime numbers $\ell$ 
such that $\ell$ is a good reduction prime for $E$ and 
$\ell \equiv 1$ (mod $p^{N}$). 
For each $\ell$, we fix a generator $\eta_{\ell}$
of $({\bf Z}/\ell {\bf Z})^{\times}$ and 
define $\log_{{\bf F}_{\ell}}(a) \in {\bf Z}/(\ell-1)$ 
by 
$\eta_{\ell}^{\log_{{\bf F}_{\ell}}(a)} \equiv a$ (mod $\ell$).  

Let $f(z)=\Sigma a_{n} e^{2 \pi i nz}$ be the modular form 
corresponding to $E$.
For a positive integer $m$ and the cyclotomic field 
${\bf Q}(\mu_{m})$, we denote by 
$\sigma_{a} \in \Gal({\bf Q}(\mu_{m})/{\bf Q})$ the 
element such that $\sigma_{a}(\zeta)=\zeta^{a}$ for any
$\zeta \in \mu_{m}$. 
We consider the modular element 
$\sum_{a=1, (a,m)=1}^{m} [\frac{a}{m}] \sigma_{a}
\in {\bf C}[\Gal({\bf Q}(\mu_{m})/{\bf Q})]$
of Mazur and Tate (\cite{MT}) 
where $[\frac{a}{m}]=2 \pi i \int_{\infty}^{a/m} f(z) dz$ is 
the usual modular symbol. 
We only consider the real part 
\begin{equation} \label{ModularElement}
\tilde{\theta}_{{\bf Q}(\mu_{m})}=
\sum_{\stackrel{\scriptstyle a=1}{(a,m)=1}}^{m} 
\frac{{\rm Re}([\frac{a}{m}])}{\Omega_{E}^{+}} \sigma_{a} 
\in {\bf Q}[\Gal({\bf Q}(\mu_{m})/{\bf Q})]
\end{equation}
where $\Omega_{E}^{+}=\int_{E({\bf R})} \omega_{E}$ is the 
N\'{e}ron period.  
Suppose that $m$ is a squarefree product of primes in ${\cal P}^{(N)}$. 
Since we are assuming the $G_{{\bf Q}}$-module $E[p]$ 
of $p$-torsion points is irreducible,  
we know $\tilde{\theta}_{{\bf Q}(\mu_{m})} \in 
{\bf Z}_{p}[\Gal({\bf Q}(\mu_{m})/{\bf Q})]$
(cf. \cite{Stevens}). 
We consider the coefficient of $\tilde{\theta}_{{\bf Q}(\mu_{m})}$ 
of ``$\prod_{\ell \vert m} (\sigma_{\eta_{\ell}}-1)$", more explicitly 
we define 
\begin{equation} \label{ModularSymbolDelta}
\tilde{\delta}_{m}=\sum_{\stackrel{\scriptstyle a=1}{(a,m)=1}}^{m} 
\frac{{\rm Re}([\frac{a}{m}])}{\Omega_{E}^{+}} 
(\prod_{\ell \mid m} \log_{{\bf F}_{\ell}}(a)) \in {\bf Z}/p^{N}
\end{equation}
where $\log_{{\bf F}_{\ell}}(a)$ means 
the image of $\log_{{\bf F}_{\ell}}(a)$ under the canonical 
homomorphism ${\bf Z}/(\ell-1) \longrightarrow 
{\bf Z}/p^{N}$. 
Let $\ord_{p}:{\bf Z}/p^{N} \longrightarrow \{0,1,...,N-1, \infty\}$ 
be the $p$-adic valuation normalized as  
$\ord_{p}(p)=1$ and $\ord_{p}(0)=\infty$. 
We note that $\ord_{p}(\tilde{\delta}_{m})$ does not 
depend on the choices of $\eta_{\ell}$ for $\ell \vert m$. 
We define $\tilde{\delta}_{1}=\theta_{{\bf Q}}
={\rm Re}([0])/\Omega_{E}^{+}
=L(E,1)/\Omega_{E}^{+}$. 

For a squarefree product $m$ of primes, 
we define $\epsilon(m)$ to be the number of prime divisors of 
$m$, namely $\epsilon(m)=r$ if $m=\ell_{1}\cdot...\cdot \ell_{r}$. 
Let ${\cal N}^{(N)}$ be the set of squarefree products of 
primes in ${\cal P}^{(N)}$. 
We suppose $1$ is in ${\cal N}^{(N)}$.
For each integer $i \geq 0$, we define the ideal 
$\Theta_{i}({\bf Q})^{(N,\delta)}$ 
of ${\bf Z}/p^{N}$ to be the ideal generated by 
all $\tilde{\delta}_{m}$ such that $\epsilon(m) \leq i$ 
for all $m \in {\cal N}^{(N)}$;
\begin{equation}
\Theta_{i}({\bf Q})^{(N,\delta)} = 
(\{\tilde{\delta}_{m} \mid \epsilon(m) \leq i \ \mbox{and} \  
m \in {\cal N}^{(N)}\}) \subset {\bf Z}/p^{N}.
\end{equation}
We define $n_{i,N}\in \{0,1,...,N-1, \infty\}$ by  
$\Theta_{i}({\bf Q})^{(N,\delta)} =p^{n_{i,N}}({\bf Z}/p^{N})$ 
(we define $n_{i,N}=\infty$ if $\Theta_{i}({\bf Q})^{(N,\delta)}=0$).

\begin{Theorem} {\rm (\cite{Ku6} Theorem B, Theorem 9.3.1 and (9.14))} 
\label{MT2}
We assume that the main conjecture for $(E, {\bf Q}_{\infty}/{\bf Q})$ 
$($see {\rm (\ref{MC})}$)$ 
and the $p$-adic height pairing is non-degenerate. \\
{\rm (1)} $n_{i,N}$ does not depend on $N$ 
when $N$ is sufficiently large (for example, when 
$N>2 \ord_{p}(\eta_{0})$ where $\eta_{0}$ is the leading term 
of the $p$-adic $L$-function, see \S 9.4 in \cite{Ku6}). 
We put $n_{i}=n_{i,N}$ for $N \gg 0$. 
In other words, we define $n_{i}$ by 
$$
{\lim\limits_{\longleftarrow}}\Theta_{i}({\bf Q})^{(N,\delta)}
=p^{n_{i}}{\bf Z}_{p} \subset {\bf Z}_{p}.$$
We denote this ideal of ${\bf Z}_{p}$ 
by $\Theta_{i}({\bf Q})^{(\delta)}$. \\
{\rm (2)}
Consider the Pontrjagin dual $\Sel(E/{\bf Q}, E[p^{\infty}])^{\vee}$ of 
the Selmer group. 
Suppose that 
$$\rank_{{\bf Z}_{p}} \Sel(E/{\bf Q}, E[p^{\infty}])^{\vee} =r  (\in {\bf Z}_{\geq 0}), \ \mbox{and} \  
\dim_{{\bf F}_{p}} \Sel(E/{\bf Q}, E[p])^{\vee} =a.$$ 
Then we have 
$$\Theta_{0}({\bf Q})^{(\delta)}=...=\Theta_{r-1}({\bf Q})^{(\delta)}
=0 \ \mbox{and} \ \Theta_{r}({\bf Q})^{(\delta)} \neq 0.$$
For any $i \geq r$, $n_{i}$ is an even number, and
$$
\begin{array}{l}
p^{n_{r}}=\#  
(\Sel(E/{\bf Q}, E[p^{\infty}])^{\vee})_{\tors},\\[3mm]
n_{a}=0, \ \mbox{and}
\end{array}
$$
$$\Sel(E/{\bf Q}, E[p^{\infty}])^{\vee} \simeq {\bf Z}_{p}^{\oplus r} 
\oplus 
({\bf Z}/p^{\frac{n_{r}-n_{r+2}}{2}})^{\oplus 2} \oplus 
({\bf Z}/p^{\frac{n_{r+2}-n_{r+4}}{2}})^{\oplus 2} \oplus
...\oplus
({\bf Z}/p^{\frac{n_{a-2}-n_{a}}{2}})^{\oplus 2}
$$
hold.
\end{Theorem}

In particular, 
knowing $\Theta_{i}({\bf Q})^{(\delta)}$ for all $i \geq 0$ 
completely determines the structure 
of $\Sel(E/{\bf Q}, E[p^{\infty}])^{\vee}$ 
as a ${\bf Z}_{p}$-module. 
Namely, the modular symbols determine the structure of 
the Selmer group under our assumptions.

\subsection{Main Results} \label{subsec12}

We define 
$${\cal P}_{1}^{(N)}=\{\ell \in {\cal P}^{(N)} \mid 
H^{0}({\bf F}_{\ell}, E[p^{N}]) 
\simeq {\bf Z}/p^{N} \}.
$$
This is an infinite set by Chebotarev density theorem 
since we are assuming (ii) (see \cite{Ku6} \S 4.3). 
We define ${\cal N}_{1}^{(N)}$ to be the set of squarefree products of 
primes in ${\cal P}_{1}^{(N)}$. 
Again, we suppose $1 \in {\cal N}_{1}^{(N)}$.
We propose the following conjecture. 

\begin{Conjecture} \label{Conj1}
There is $m \in {\cal N}_{1}^{(N)}$ such that 
$\tilde{\delta}_{m}$ is a unit in ${\bf Z}/p^{N}$, namely 
$$\ord_{p}(\tilde{\delta}_{m})=0.$$
\end{Conjecture}

Numerically, it is easy to compute $\tilde{\delta}_{m}$, 
so it is easy to check this conjecture. 

\begin{Theorem} {\rm (\cite{Ku6} Theorem 9.3.1)}
If we assume the main conjecture and 
the non-degeneracy of the $p$-adic height pairing, 
Conjecture \ref{Conj1} holds true. 
\end{Theorem}

In fact, we obtain Conjecture \ref{Conj1}, 
considering the case $i=a$ in Theorem \ref{MT2} 
(cf. $i=s$ in Theorem 9.3.1 in \cite{Ku6}).

\vspace{5mm}

{\it From now on, we do not assume the main conjecture (\ref{MC}) 
nor the non-degeneracy of the $p$-adic height pairing}. 

We define the Selmer group $\Sel({\bf Z}[1/m], E[p^{N}])$ 
by 
$$\Sel({\bf Z}[1/m], E[p^{N}]) = 
\Ker(H^{1}({\bf Q}, E[p^{N}]) \longrightarrow 
\prod_{v \not \hspace{0.7mm} \mid m} H^{1}({\bf Q}_{v}, E[p^{N}])
/E({\bf Q}_{v}) \otimes {\bf Z}/p^{N}).$$
If all bad primes and $p$ divide $m$, we know 
$\Sel({\bf Z}[1/m], E[p^{N}])$ is equal to the \'{e}tale cohomology 
group $H^{1}_{et}(\Spec {\bf Z}[1/m],E[p^{N}])$, which 
explains the notation ``$\Sel({\bf Z}[1/m], E[p^{N}])$". 
(We use $\Sel({\bf Z}[1/m], E[p^{N}])$ for $m \in {\cal N}_{1}^{(N)}$ 
in this paper, but 
$E[p^{N}]$ is not an \'{e}tale sheaf on $\Spec {\bf Z}[1/m]$ 
for such $m$.)  

Let $\lambda$ be the $\lambda$-invariant of 
$\Sel(E/{\bf Q}_{\infty}, E[p^{\infty}])^{\vee}$. 
We put $n_{\lambda} = \min\{n \in {\bf Z} \mid p^{n}-1 \geq \lambda\}$ 
and $d_{n}=n_{\lambda}+Nn$ for $n \in {\bf Z}_{\geq 0}$. 
We define 
\begin{equation} \label{PNn}
{\cal P}_{1}^{(N,n)} = \{ \ell \in {\cal P}_{1}^{(N)}
 \mid \ell \equiv 1 \ \mbox{(mod $p^{d_{n}}$)} \}
\end{equation}
(then 
${\cal P}_{1}^{(N,n)} \subset 
{\cal P}_{1}^{(N)}({\bf Q}_{[n]})$ holds, 
see the end of \S \ref{subsec31} for this fact, and see \S \ref{subsec31} for 
the definition of the set ${\cal P}_{1}^{(N)}({\bf Q}_{[n]})$). 
We denote by ${\cal N}_{1}^{(N,n)}$ the set of squarefree products of 
primes in ${\cal P}_{1}^{(N,n)}$. 

In this paper, 
for any finite abelian $p$-extension $K/{\bf Q}$ in which 
all bad primes of $E$ are unramified, 
we prove in \S \ref{section4} the following theorem for 
${\bf Z}/p^{N}[\Gal(K/{\bf Q})]$-modules 
$\Sel(E/K, E[p^{N}])$ and 
$\Sel(O_{K}[1/m], E[p^{N}])$
(see Corollary \ref{injectivitytheorem}
and Theorem \ref{RelationTheorem}). 
We simply state it in the case $K={\bf Q}$ below. 
An essential ingredient in this paper is the Kolyvagin system of 
Gauss sum type. 
We construct Kolyvagin systems $\kappa_{m,\ell} 
\in \Sel({\bf Z}[1/m\ell], E[p^{N}])$ 
for $(m,\ell)$ satisfying 
$\ell \in {\cal P}_{1}^{(N,\epsilon(m\ell)+1)}$
and 
$m \ell \in {\cal N}_{1}^{(N,\epsilon(m\ell)+1)}$
(see \S \ref{subsec34} and Propositions \ref{KKP2}) 
by the method in \cite{Ku6}.
(We can construct these elements, using 
the half of the main conjecture proved by Kato \cite{KK1}.)
The essential difference between our 
Kolyvagin systems $\kappa_{m, \ell}$ of Gauss sum type 
and Kolyvagin systems in Mazur and Rubin \cite{MR}
is that our $\kappa_{m, \ell}$ is related to $L$-values. 
In particular, $\kappa_{m, \ell}$ satisfies a remarkable 
property $\phi_{\ell}(\kappa_{m, \ell})=-\delta_{m\ell}t_{\ell,K}$ 
(see Propositions \ref{KKP2} (4)) though we do not 
explain the notation here. 

\begin{Theorem} \label{MMT1}
Assume that $\ord_{p}(\tilde{\delta}_{m})=0$ for some 
$m \in {\cal N}_{1}^{(N)}$. \\
{\rm (1)} The canonical homomorphism 
$$s_{m}: \Sel(E/{\bf Q}, E[p^{N}]) \longrightarrow 
\bigoplus_{\ell \mid m} E({\bf Q}_{\ell})\otimes{\bf Z}/p^{N} 
\simeq 
\bigoplus_{\ell \mid m} E({\bf Q}_{\ell})\otimes{\bf Z}/p^{N} 
\simeq 
({\bf Z}/p^{N})^{\epsilon(m)}$$
is injective. \\
{\rm (2)}
Assume further that $m \in {\cal N}_{1}^{(N,\epsilon(m)+1)}$ and 
that $m$ is admissible (for the definition of the notion ``admissible", see
the paragraph before Proposition \ref{KKP1}).
Then $\Sel({\bf Z}[1/m], E[p^{N}])$ is a free ${\bf Z}/p^{N}$-module 
of rank $\epsilon(m)$, and 
$\{\kappa_{\frac{m}{\ell},\ell}\}_{\ell \mid m}$ is a basis of 
$\Sel({\bf Z}[1/m], E[p^{N}])$. \\
{\rm (3)} 
We define a matrix
${\cal A}$ as in {\rm (\ref{relationmatrix})} 
in Theorem {\rm \ref{RelationTheorem}}, using 
$\kappa_{\frac{m}{\ell},\ell}$. 
Then ${\cal A}$ is a 
relation matrix of the Pontrjagin dual $\Sel(E/{\bf Q}, E[p^{N}])^{\vee}$
of the Selmer group; 
namely if $f_{\cal A}: ({\bf Z}/p^{N})^{\epsilon(m)} 
\longrightarrow ({\bf Z}/p^{N})^{\epsilon(m)}$
is the homomorphism corresponding to 
the above matrix ${\cal A}$, then we have 
$$\Coker (f_{\cal A}) \simeq \Sel(E/{\bf Q}, E[p^{N}])^{\vee}.$$
\end{Theorem}

It is worth noting that we get nontrivial 
(moreover, linearly independent) elements 
in the Selmer groups.

\vspace{5mm}

The ideals $\Theta_{i}({\bf Q})^{(\delta)}$ in 
Theorem \ref{MT2} are not suitable for numerical computations 
because we have to compute {\it infinitely many} $\tilde{\delta}_{m}$. 
On the other hand, we can easily find $m$ with 
$\ord_{p}(\tilde{\delta}_{m})=0$ numerically. 
Since $s_{m}$ is injective, we can get information of 
the Selmer group from the image of $s_{m}$, which 
is an advantage of Theorem \ref{MMT1} and the next Theorem \ref{MMT2} 
(see also the comment in the 
end of Example (5) in \S \ref{subsection53}). 

\vspace{5mm}

We next consider the case $N=1$, so $\Sel(E/{\bf Q}, E[p])$. 
Now we regard $\tilde{\delta}_{m}$ as an element of  ${\bf F}_{p}$ 
for $m \in {\cal N}_{1}^{(1)}$. 
We say $m$ is {\it $\delta$-minimal} if 
$\tilde{\delta}_{m} \neq 0$ and $\tilde{\delta}_{d}=0$ 
for all divisors $d$ of $m$ with $1 \leq d<m$. 
Our next conjecture claims that the structure (the dimension) of 
$\Sel(E/{\bf Q}, E[p])$
is determined by a $\delta$-minimal $m$, 
therefore can be easily computed numerically. 

\begin{Conjecture} \label{Conj2}
If $m \in {\cal N}_{1}^{(1)}$ is $\delta$-minimal,  
the canonical homomorphism 
$$s_{m}: \Sel(E/{\bf Q}, E[p]) \longrightarrow 
\bigoplus_{\ell \mid m} E({\bf Q}_{\ell})\otimes{\bf Z}/p 
\simeq 
\bigoplus_{\ell \mid m} E({\bf F}_{\ell})\otimes{\bf Z}/p 
\simeq 
({\bf Z}/p^{N})^{\epsilon(m)}$$
is bijective. 
In particular, $\dim_{{\bf F}_{p}}\Sel(E/{\bf Q}, E[p])=\epsilon(m)$. 
\end{Conjecture}
 
If $m \in {\cal N}_{1}^{(1)}$ is $\delta$-minimal, 
the above homomorphism $s_{m}: \Sel(E/{\bf Q}, E[p]) \longrightarrow 
({\bf Z}/p^{N})^{\epsilon(m)}$ 
is injective by Theorem \ref{MMT1} (1), 
so we know  
$$\dim_{{\bf F}_{p}}\Sel(E/{\bf Q}, E[p]) \leq \epsilon(m).$$
Therefore, the problem is in showing the other inequality. 

We note that the analogue of the above conjecture for ideal class groups 
does not hold (see \S \ref{subsection54}). 
But we hope that Conjecture \ref{Conj2} holds for the Selmer groups of 
elliptic curves.
We construct in \S \ref{section5} a modified version 
$\kappa_{m,\ell}^{q,q',z}$ of Kolyvagin systems of Gauss sum type  
for any $(m,\ell)$ with $m \ell \in {\cal N}_{1}^{(N)}$. 
(The Kolyvagin system $\kappa_{m,\ell}$ in \S 3 is defined for $(m,\ell)$ with 
$m \ell \in {\cal N}_{1}^{(N,\epsilon(m\ell)+1)}$, 
but $\kappa_{m,\ell}^{q,q',z}$ is defined for more general 
$(m,\ell)$, namely for $(m,\ell)$ with $m \ell \in {\cal N}_{1}^{(N)}$.)
Using the modified Kolyvagin 
system $\kappa_{m,\ell}^{q,q',z}$, 
we prove the following. 

\begin{Theorem} \label{MMT2}
{\rm (1)} If $\epsilon(m)=0$, $1$, then Conjecture {\rm \ref{Conj2}} 
is true. \\
{\rm (2)} If there is $\ell \in  {\cal P}^{(1)}$ which is $\delta$-minimal
(so $\epsilon(\ell)=1$), then 
$$\Sel(E/{\bf Q}, E[p^{\infty}]) \simeq {\bf Q}_{p}/{\bf Z}_{p}.$$ 
Moreover, if there is $\ell \in {\cal P}_{1}^{(1)}$ 
which is $\delta$-minimal and which satisfies 
$\ell \equiv 1$ (mod  $p^{n_{\lambda'}+2})$ 
where $\lambda'$ is the analytic $\lambda$-invariant of 
$(E, {\bf Q}_{\infty}/{\bf Q})$, 
then the main conjecture {\rm (\ref{MC})} for 
$\Sel(E/{\bf Q}_{\infty}, E[p^{\infty}])$ 
holds true. 
In this case, $\Sel(E/{\bf Q}_{\infty}, E[p^{\infty}])^{\vee}$ 
is generated by one element as a 
${\bf Z}_{p}[[\Gal({\bf Q}_{\infty}/{\bf Q})]]$-module. \\
{\rm (3)} If $\epsilon(m)=2$ and $m$ is admissible, 
then Conjecture {\rm \ref{Conj2}} is true. \\
{\rm (4)} Suppose that $\epsilon(m)=3$ and $m=\ell_{1}\ell_{2}\ell_{3}$. 
Assume that $m$ is admissible and the natural maps 
$s_{\ell_{i}}: \Sel(E/{\bf Q}, E[p]) \longrightarrow 
E({\bf F}_{\ell_{i}})\otimes{\bf Z}/p$ are surjective both for 
$i=1$ and $i=2$. 
Then Conjecture {\rm \ref{Conj2}} is true.
\end{Theorem}

In this way, we can determine the Selmer groups by finite numbers of computations 
in several cases. 
We give several numerical examples in \S \ref{subsection52}.

\begin{Remark}
\begin{rm}
Concerning the Fitting ideals and the annihilator ideals 
of some Selmer groups, 
we prove the following in this paper. 
Let $K/{\bf Q}$ be a finite abelian 
$p$-extension in which all bad primes of $E$ are unramified.
We take a finite set $S$ of good reduction primes, which contains 
all ramifying primes in $K/{\bf Q}$ except $p$. 
Let $m$ be the product of primes in $S$. 
We prove that the initial 
Fitting ideal of the $R_{K}={\bf Z}_{p}[\Gal(K/{\bf Q})]$-module 
$\Sel(O_{K}[1/m], E[p^{\infty}])^{\vee}$ is principal, 
and 
$$
\xi_{K,S} \in  \Fitt_{0, R_{K}}
(\Sel(O_{K}[1/m], E[p^{\infty}])^{\vee})
$$
where $\xi_{K,S}$ is an element of $R_{K}$ which is 
explicitly constructed from modular symbols 
(see (\ref{FittingIdeal02})). 
If the main conjecture (\ref{MC}) for 
$(E, {\bf Q}_{\infty}/{\bf Q})$ holds, the equality 
$\Fitt_{0, R_{K}}
(\Sel(O_{K}[1/m], E[p^{\infty}])^{\vee})=\xi_{K,S}R_{K}$
holds (see Remark \ref{FittingIdeal00}). 
We prove the Iwasawa theoretical version in 
Theorem \ref{T1}.

Let $\vartheta_{K}$ be the image of the $p$-adic $L$-function, 
which is also explicitly constructed from modular symbols. 
We show in Theorem \ref{AnnTh1}
$$\vartheta_{K} \in  \Ann_{R_{K}}
(\Sel(O_{K}[1/m], E[p^{\infty}])^{\vee}).$$

Concerning the higher Fitting ideals (cf. \S \ref{subsec24}), 
we show 
$$\tilde{\delta}_{m} \in \Fitt_{\epsilon(m), {\bf Z}/p^{N}}
(\Sel(E/{\bf Q}, E[p^{N}])^{\vee})$$ 
where $\Fitt_{i, R}(M)$ is the $i$-th Fitting ideal of an $R$-module 
$M$. 
We prove a slightly generalized version for 
$K$ which is in the cyclotomic ${\bf Z}_{p}$-extension 
${\bf Q}_{\infty}$ of ${\bf Q}$ 
(see Theorem \ref{HigherFittingIdeal} and 
Corollary \ref{HigherFittingIdeal2}).
\end{rm}
\end{Remark}

\vspace{5mm}

I would like to thank John Coates heartily for his helpful advice and 
for discussion with him, especially for the discussion in March 2013,  
which played an essential role in my producing this paper. 
I also thank heartily Kazuya Kato for 
his constant interest in the results of this paper.
I also thank Kazuo Matsuno and Christian Wuthrich very much 
for their helping me to compute modular symbols.

\section{Selmer groups and $p$-adic $L$-functions} \label{section2}

\subsection{Modular symbols and $p$-adic $L$-functions} \label{subsec21}

Let $E$ be an elliptic curve over ${\bf Q}$, and 
$f(z)=\Sigma a_{n} e^{2 \pi i nz}$ the modular form 
corresponding to $E$. 
In this section, we assume that $p$ is 
a prime number satisfying (i), (ii), (iii) in 
\S \ref{subsec11}.
We define ${\cal P}_{good}=\{\ell \mid $ $\ell$ is a 
good reduction prime for $E$ $\} \setminus \{p\}$.  
For any finite abelian extension $K/{\bf Q}$, 
we denote by $K_{\infty}/K$ the cyclotomic ${\bf Z}_{p}$-extension. 
For a real abelian field $K$ of conductor $m$, 
we define $\tilde{\theta}_{K}$ to be the image of 
$\tilde{\theta}_{{\bf Q}(\mu_{m})}$ in ${\bf Q}[\Gal(K/{\bf Q})]$ 
where $\tilde{\theta}_{{\bf Q}(\mu_{m})}$ is defined in 
(\ref{ModularElement}).

We write 
$$R_{K}={\bf Z}_{p}[\Gal(K/{\bf Q})]  \ \mbox{and} \
\Lambda_{K_{\infty}}={\bf Z}_{p}[[\Gal(K_{\infty}/{\bf Q})]].$$ 

For any positive integer $n$, 
we simply write 
$R_{{\bf Q}(\mu_{n})}=R_{n}$ in this subsection.
For any positive integers $d$, $c$ such that $d \mid c$, 
we define the norm map $\nu_{c,d}: 
R_{d}={\bf Z}_{p}[\Gal({\bf Q}(\mu_{d})/{\bf Q})] \longrightarrow 
R_{c}={\bf Z}_{p}[\Gal({\bf Q}(\mu_{c})/{\bf Q})]$ 
by 
$\sigma \mapsto \sum \tau$ where for $\sigma 
\in \Gal({\bf Q}(\mu_{d})/{\bf Q})$, 
$\tau$ runs over all elements of 
$\Gal({\bf Q}(\mu_{c})/{\bf Q})$ such that the restriction of 
$\tau$ to ${\bf Q}(\mu_{d})$ is $\sigma$.
Let $m$ be a squarefree product of primes in 
${\cal P}_{good}$, and 
$n$ a positive integer. 
By our assumption (ii), 
we know $\tilde{\theta}_{{\bf Q}(\mu_{mp^{n}})} \in 
R_{mp^{n}}$
(cf. \cite{Stevens}). 
Let $\alpha \in {\bf Z}_{p}^{\times}$ be the unit root of 
$x^{2}-a_{p}x+p=0$ and put 
$$\vartheta_{{\bf Q}(\mu_{mp^{n}})}=
\alpha^{-n}(\tilde{\theta}_{{\bf Q}(\mu_{mp^{n}})}-
\alpha^{-1}\nu_{mp^{n},mp^{n-1}}
(\tilde{\theta}_{{\bf Q}(\mu_{mp^{n-1}})})) 
\in R_{mp^{n}}$$
as usual.
Then $\{\vartheta_{{\bf Q}(\mu_{mp^{n}})}\}_{n \geq 1}$ 
is a projective system 
(cf. Mazur and Tate \cite{MT} the equation (4) on page 717) 
and we obtain an element 
$\vartheta_{{\bf Q}(\mu_{mp^{\infty}})}
\in \Lambda_{{\bf Q}(\mu_{mp^{\infty}})}$, 
which is the $p$-adic $L$-function of Mazur and Swinnerton-Dyer. 

We also use the notation $\Lambda_{np^{\infty}}=
\Lambda_{{\bf Q}(\mu_{np^{\infty}})}$ for simplicity. 
Suppose that a prime $\ell$ does not divide $mp$, and 
$c_{m\ell, m}:  \Lambda_{m\ell p^{\infty}} 
\longrightarrow \Lambda_{mp^{\infty}}$ 
is the natural restriction map. 
Then we know 
\begin{equation} \label{111}
c_{m\ell, m}(\vartheta_{{\bf Q}(\mu_{m\ell p^{\infty}})})
= (a_{\ell}-\sigma_{\ell}-\sigma_{\ell}^{-1})
\vartheta_{{\bf Q}(\mu_{mp^{\infty}})}
\end{equation}
(cf. Mazur and Tate \cite{MT} the equation (1) on page 717). 

We will construct a slightly modified element 
$\xi_{{\bf Q}(\mu_{mp^{\infty}})}$ in 
$\Lambda_{mp^{\infty}}$. 
We put $P'_{\ell}(x)=x^{2}-a_{\ell}x+\ell$. 
Let $m$ be a squarefree product of ${\cal P}_{good}$. 
For any divisor $d$ of $m$ and a prime divisor $\ell$ of $m/d$, 
$\sigma_{\ell} \in \Gal({\bf Q}(\mu_{dp^{\infty}})/{\bf Q})
= {\lim\limits_{\longleftarrow}} \Gal({\bf Q}(\mu_{dp^{n}})/{\bf Q})$ 
is defined as the projective limit of $\sigma_{\ell} 
\in \Gal({\bf Q}(\mu_{dp^{n}})/{\bf Q})$.  
We consider $P'_{\ell}(\sigma_{\ell}) \in \Lambda_{dp^{\infty}}$.
Note that 
\begin{equation} \label{112}
-\sigma_{\ell}^{-1}=(-\sigma_{\ell}^{-1}P_{\ell}'(\sigma_{\ell})
-(a_{\ell}-\sigma_{\ell}-\sigma_{\ell}^{-1}))/(\ell-1) 
\in \Lambda_{dp^{\infty}}.
\end{equation}
We put $\alpha_{d,m}=(\prod_{\ell \vert \frac{m}{d}} (-\sigma_{\ell}^{-1})) 
\vartheta_{{\bf Q}(\mu_{d p^{\infty}})}
\in \Lambda_{dp^{\infty}}$ and 
$$\xi_{{\bf Q}(\mu_{m p^{\infty}})}
=\sum_{d \vert m} \nu_{m,d}(\alpha_{d,m}) \in 
\Lambda_{mp^{\infty}}$$
where $\nu_{m,d}:\Lambda_{dp^{\infty}} 
\longrightarrow \Lambda_{mp^{\infty}}$ 
is the norm map defined similarly as above. 
(This modification $\xi_{{\bf Q}(\mu_{m p^{\infty}})}$ 
is done by the same spirit as Greither \cite{Grei3} in 
which the Deligne-Ribet $p$-adic $L$-functions are 
treated.) 
Suppose that $\ell \in {\cal P}_{good}$ is prime to $m$. 
Then by the definition of $\xi_{{\bf Q}(\mu_{mp^{\infty}})}$ 
and (\ref{111}) and (\ref{112}), we have  
\begin{eqnarray} \label{Res2}
c_{m\ell, m}(\xi_{{\bf Q}(\mu_{m\ell p^{\infty}})})
&=&
c_{m\ell, m}
(\sum_{d \vert m} 
\nu_{m\ell, d}(\alpha_{d,m \ell})
+\sum_{d \vert m} 
\nu_{m\ell, d \ell}(\alpha_{d\ell,m \ell})) \nonumber \\
&=&(\ell-1)\sum_{d \vert m} \nu_{m,d}
(-\sigma_{\ell}^{-1}\alpha_{d,m})
+\sum_{d \vert m}\nu_{m,d}(c_{d\ell, d}
(\alpha_{d\ell,m \ell})) \nonumber \\
&=&(\ell-1)\sum_{d \vert m} \nu_{m,d}
(-\sigma_{\ell}^{-1}\alpha_{d,m})
+\sum_{d \vert m}\nu_{m,d}(
(a_{\ell}-\sigma_{\ell}-\sigma_{\ell}^{-1})
\alpha_{d,m}) \nonumber \\
&=&
(-\sigma_{\ell}^{-1}P_{\ell}'(\sigma_{\ell}))
\sum_{d \vert m} \nu_{m,d}(\alpha_{d,m}) \nonumber \\
&=&
(-\sigma_{\ell}^{-1}P_{\ell}'(\sigma_{\ell}))
\xi_{{\bf Q}(\mu_{m p^{\infty}})}.
\end{eqnarray}

We denote by $\vartheta_{{\bf Q}(\mu_{m})} \in R_{{\bf Q}(\mu_{m})}$
the image of 
$\vartheta_{{\bf Q}(\mu_{mp^{\infty}})}$ under the natural map 
$\Lambda_{{\bf Q}(\mu_{mp^{\infty}})} 
\longrightarrow 
R_{{\bf Q}(\mu_{m})}$. 
We have 
\begin{equation} \label{Res3}
\vartheta_{{\bf Q}(\mu_{m})}=(1-\frac{\sigma_{p}}{\alpha})
(1-\frac{\sigma_{p}^{-1}}{\alpha})\tilde{\theta}_{{\bf Q}(\mu_{m})}.
\end{equation}
Since we are assuming $a_{p} \not \equiv 1$ (mod $p$), 
we also have $\alpha \not \equiv 1$ (mod $p$), so 
$(1-\frac{\sigma_{p}}{\alpha})
(1-\frac{\sigma_{p}^{-1}}{\alpha})$ is a unit in 
$R_{{\bf Q}(m)}$ where ${\bf Q}(m)$ is 
the maximal $p$-subextension of ${\bf Q}$ in ${\bf Q}(\mu_{m})$.

\subsection{Selmer groups} \label{subsec22}

For any algebraic extension $F/{\bf Q}$, we denote by 
$O_{F}$ the integral closure of ${\bf Z}$ in $F$.
For a positive integer $m>0$, we define a Selmer group 
$\Sel(O_{F}[1/m], E[p^{\infty}])$ by 
$$\Sel(O_{F}[1/m], E[p^{\infty}])=
\Ker(H^{1}(F, E[p^{\infty}]) \longrightarrow 
\prod_{v \not \hspace{0.7mm} \mid m} H^{1}(F_{v}, E[p^{\infty}])
/E(F_{v}) \otimes {\bf Q}_{p}/{\bf Z}_{p})$$
where $v$ runs over all primes of $F$ which are prime to $m$.  
Similarly, for a positive integer $N$, we define 
$\Sel(O_{F}[1/m], E[p^{N}])$ by 
$$\Sel(O_{F}[1/m], E[p^{N}])=
\Ker(H^{1}(F, E[p^{N}]) \longrightarrow 
\prod_{v \not \hspace{0.7mm} \mid m} H^{1}(F_{v}, E[p^{N}])
/E(F_{v}) \otimes {\bf Z}/p^{N}).$$
In the case $m=1$, we denote them by 
$\Sel(O_{F}, E[p^{\infty}])$, $\Sel(O_{F}, E[p^{N}])$, 
which are classical Selmer groups. 
We also use the notation $\Sel(E/F, E[p^{\infty}])$, 
$\Sel(E/F, E[p^{N}])$ for them, namely 
$$\Sel(E/F, E[p^{\infty}])=\Sel(O_{F}, E[p^{\infty}]), \ \ 
\Sel(E/F, E[p^{N}])=\Sel(O_{F}, E[p^{N}]).$$ 

For a finite abelian extension $K/{\bf Q}$, we denote by 
$K_{\infty}/K$ the cyclotomic ${\bf Z}_{p}$-extension, and 
put $\Lambda_{K_{\infty}}={\bf Z}_{p}[[\Gal(K_{\infty}/{\bf Q})]]$.
The Pontrjagin dual 
$\Sel(O_{K_{\infty}}, E[p^{\infty}])^{\vee}$
is a torsion $\Lambda_{K_{\infty}}$-module 
(Kato \cite{KK1} Theorem 17.4). 

When the conductor of $K$ is $m$, 
we define $\vartheta_{K_{\infty}} \in \Lambda_{K_{\infty}}$ 
to be the image of $\vartheta_{{\bf Q}(\mu_{mp^{\infty}})}$, 
and also $\xi_{K_{\infty}} \in \Lambda_{K_{\infty}}$ 
to be the image of $\xi_{{\bf Q}(\mu_{mp^{\infty}})}$.
The Iwasawa main conjecture for $(E, {\bf Q}_{\infty}/{\bf Q})$ 
is the equality between the characteristic ideal of the Selmer 
group and the ideal generated by the $p$-adic $L$-function; 
\begin{equation} \label{MC}
\cc(\Sel(O_{{\bf Q}_{\infty}}, E[p^{\infty}])^{\vee})=
\vartheta_{{\bf Q}_{\infty}}\Lambda_{{\bf Q}_{\infty}}.
\end{equation}
Since we are assuming the Galois action on the Tate module is 
surjective, 
we know 
$\vartheta_{{\bf Q}_{\infty}} \in 
\cc(\Sel(O_{{\bf Q}_{\infty}}, E[p^{\infty}])^{\vee})$
by Kato \cite{KK1} Theorem 17.4. 
Skinner and Urban \cite{SU1} proved the equality (\ref{MC}) 
under mild conditions. 
Namely, 
under our assumptions (i), (ii), 
they proved the main conjecture (\ref{MC}) if 
there is a bad prime $\ell$ which is ramified in ${\bf Q}(E[p])$ 
(\cite{SU1} Theorem 3.33).

More generally, let $\psi$ be an even Dirichlet character 
and $K$ be the abelian field corresponding to the kernel of 
$\psi$, namely $K$ is the field such that $\psi$ induces 
a faithful character of $\Gal(K/{\bf Q})$.  
We assume $K \cap {\bf Q}_{\infty}={\bf Q}$.
In this paper, for any finite abelian $p$-group $G$, 
any ${\bf Z}_{p}[G]$-module $M$ and 
any character $\psi:G \longrightarrow \overline 
{{\bf Q}}_{p}^{\times}$, 
we define the $\psi$-quotient $M_{\psi}$ by 
$M \otimes_{{\bf Z}_{p}[G]} O_{\psi}$
where $O_{\psi}={\bf Z}_{p}[\ii \psi]$ which is regarded as 
a ${\bf Z}_{p}[G]$-module by 
$\sigma x=\psi(\sigma)x$ for any $\sigma \in G$ 
and $x \in O_{\psi}$. 
We consider  
$(\Sel(O_{K_{\infty}}, E[p^{\infty}])^{\vee})_{\psi}$, 
which is a $\Lambda_{\psi}$-module where
$\Lambda_{\psi}=(\Lambda_{K_{\infty}})_{\psi}
= O_{\psi}[[\Gal(K_{\infty}/K)]]$.   
We denote the image of $\vartheta_{K_{\infty}}$ in 
$\Lambda_{\psi}$ by $\psi(\vartheta_{K_{\infty}})$. 
Then the main conjecture states
\begin{equation} \label{MC2}
\cc((\Sel(O_{K_{\infty}}, E[p^{\infty}])^{\vee})_{\psi})=
\psi(\vartheta_{K_{\infty}})\Lambda_{\psi}.
\end{equation}
We also note that 
$\psi(\vartheta_{K_{\infty}})\Lambda_{\psi}=
\psi(\xi_{K_{\infty}})\Lambda_{\psi}$.
By Kato \cite{KK1}, we know $\psi(\vartheta_{K_{\infty}})$, 
$\psi(\xi_{K_{\infty}}) \in 
\cc((\Sel(O_{K_{\infty}}, E[p^{\infty}])^{\vee})_{\psi})$.

Let $S \subset {\cal P}_{good}$ be a finite set of good primes, 
and $K/{\bf Q}$ be a finite abelian extension. 
We denote by $S_{{\rm ram}}(K)$ the subset of $S$ which 
consists of all ramifying primes in $K$ inside $S$. 
Recall that $P'_{\ell}(x)=x^{2}-a_{\ell}x+\ell$. 
We define 
$$\xi_{K_{\infty},S}=\xi_{K_{\infty}} 
\prod_{\ell \in S \setminus S_{{\rm ram}}(K)}
(-\sigma_{\ell}^{-1}P_{\ell}'(\sigma_{\ell})).$$
So $\xi_{K_{\infty},S}=\xi_{K_{\infty}}$ 
if $S$ contains only ramifying primes in $K$.
Suppose that $S$ contains all ramifying primes in $K$ 
and $F$ is a subfield of $K$. 
We denote by $c_{K_{\infty}/F_{\infty}}:\Lambda_{K_{\infty}} 
\longrightarrow \Lambda_{F_{\infty}}$ the natural 
restriction map. 
Using (\ref{Res2}) and the above 
definition of $\xi_{K_{\infty},S}$, 
we have 
\begin{equation} \label{Res4}
c_{K_{\infty}/F_{\infty}}(\xi_{K_{\infty},S})=\xi_{F_{\infty},S}.
\end{equation}
 
For any positive integer $m$ whose prime divisors are in ${\cal P}_{good}$, 
we have an exact sequence 
$$0 \longrightarrow  \Sel(O_{K_{\infty}}, E[p^{\infty}])
\longrightarrow \Sel(O_{K_{\infty}}[1/m], E[p^{\infty}])
\longrightarrow  \bigoplus_{v \mid m} H^{1}(K_{\infty,v}, E[p^{\infty}])
\longrightarrow  0$$
because $E(K_{\infty,v})\otimes{\bf Q}_{p}/{\bf Z}_{p}=0$ 
(for the surjectivity of the third map, see Greenberg Lemma 4.6 
in \cite{Gr2}). 
For a prime $v$ of $K_{\infty}$, let $K_{\infty,v,nr}/K_{\infty,v}$ be
the maximal unramified extension, and $\Gamma_{v}=
\Gal(K_{\infty,v,nr}/K_{\infty,v})$. 
Suppose $v$ divides $m$. 
Since $v$ is a good reduction prime, we have 
$H^{1}(K_{\infty,v}, E[p^{\infty}])=
\Hom_{Cont}(G_{K_{\infty,v,nr}}, E[p^{\infty}])^{\Gamma_{v}}
=E[p^{\infty}](-1)^{\Gamma_{v}}$ where $(-1)$ is the Tate twist. 
By the Weil pairing, the Pontrjagin dual of $E[p^{\infty}](-1)$ 
is the Tate module $T_{p}(E)$. 
Therefore, taking the Pontrjagin dual of the above exact 
sequence, we have an exact sequence 
\begin{equation} \label{ES1}
0 \longrightarrow  \bigoplus_{v \mid m} T_{p}(E)_{\Gamma_{v}} 
\longrightarrow \Sel(O_{K_{\infty}}[1/m], E[p^{\infty}])^{\vee}
\longrightarrow  \Sel(O_{K_{\infty}}, E[p^{\infty}])^{\vee}
\longrightarrow  0.
\end{equation}
Note that $T_{p}(E)_{\Gamma_{v}}$ is free over ${\bf Z}_{p}$ 
because $\Gamma_{v}$ is profinite of order prime to $p$. 

\vspace{5mm}

Let $K/{\bf Q}$ be a finite abelian 
$p$-extension in which all bad primes of $E$ are unramified.
Suppose that $S$ is a finite subset of ${\cal P}_{good}$ such that
$S$ contains all ramifying primes in $K/{\bf Q}$ except $p$. 
Let $m$ be a squarefree product of all primes in $S$.  

\begin{Theorem} \label{GreenbergAMS} {\rm (Greenberg)}
$\Sel(O_{K_{\infty}}[1/m], E[p^{\infty}])^{\vee}$ is 
of projective dimension $\leq 1$ as a $\Lambda_{K_{\infty}}$-module.
\end{Theorem}

This is proved by Greenberg in \cite{Gr3} Theorem 1 
(the condition (iv) in \S \ref{subsec11} in this paper 
is not needed here, 
see also Proposition 3.3.1 in \cite{Gr3}). 
For more general $p$-adic representations, this is proved in \cite{Ku6} Proposition 1.6.7.
We will give a sketch of the proof because some results in the proof 
will be used later. 

Since we can take some finite abelian extension $K'/{\bf Q}$ such that 
$K_{\infty}=K'_{\infty}$ and $K' \cap {\bf Q}_{\infty} = {\bf Q}$, 
we may assume that $K \cap {\bf Q}_{\infty} = {\bf Q}$
and $p$ is unramified in $K$. 
Since we are assuming that $E[p]$ is an irreducible $G_{{\bf Q}}$-module, 
we know that $\Sel(O_{K_{\infty}}, E[p^{\infty}])^{\vee}$ 
has no nontrivial finite ${\bf Z}_{p}[[\Gal(K_{\infty}/K)]]$-submodule 
by Greenberg (\cite{Gr2} Propositions 4.14, 4.15). 
We also assumed that the $\mu$-invariant of 
$\Sel(O_{{\bf Q}_{\infty}}, E[p^{\infty}])^{\vee}$
is zero, which implies the vanishing of the $\mu$-invariant of 
$\Sel(O_{K_{\infty}}, E[p^{\infty}])^{\vee}$
by Hachimori and Matsuno \cite{HM1}.  
Therefore, $\Sel(O_{K_{\infty}}, E[p^{\infty}])^{\vee}$ 
is a free ${\bf Z}_{p}$-module of finite rank. 
By the exact sequence (\ref{ES1}), 
$\Sel(O_{K_{\infty}}[1/m], E[p^{\infty}])^{\vee}$ is also 
a free ${\bf Z}_{p}$-module of finite rank. 

Put $G=\Gal(K/{\bf Q})$. 
By the definition of the Selmer group and our assumption that 
all primes dividing $m$ are good reduction primes, 
we have 
$\Sel(O_{K_{\infty}}[1/m], E[p^{\infty}])^{G}=
\Sel(O_{{\bf Q}_{\infty}}[1/m], E[p^{\infty}])$.
Since we assumed that the $\mu$-invariant is zero, 
$\Sel(O_{{\bf Q}_{\infty}}[1/m], E[p^{\infty}])$ is divisible. 
This shows that the corestriction map 
$\Sel(O_{K_{\infty}}[1/m], E[p^{\infty}]) \longrightarrow 
\Sel(O_{{\bf Q}_{\infty}}[1/m], E[p^{\infty}])$ is 
surjective. 
Therefore, ${\hat H}^{0}(G, 
\Sel(O_{K_{\infty}}[1/m], E[p^{\infty}]))=0$. 

Next we will show that $H^{1}(G, 
\Sel(O_{K_{\infty}}[1/m], E[p^{\infty}]))=0$.
Let $N_{E}$ be the conductor of $E$ and 
put $m'=mpN_{E}$. 
We know 
$\Sel(O_{K_{\infty}}[1/m'], E[p^{\infty}])$ is equal to 
the \'{e}tale cohomology group 
$H^{1}_{et}(\Spec O_{K_{\infty}}[1/m'],E[p^{\infty}])$. 
We have an exact sequence
\begin{equation} \label{ES2}
0 \longrightarrow \Sel(O_{K_{\infty}}[1/m], E[p^{\infty}])
\longrightarrow  \Sel(O_{K_{\infty}}[1/m'], E[p^{\infty}])
\longrightarrow \bigoplus_{v \mid \frac{m'}{m}} H^{2}_{v}(K_{\infty,v})
\longrightarrow  0
\end{equation}
where $H^{2}_{v}(K_{\infty,v})=H^{1}(K_{\infty,v}, E[p^{\infty}])/
(E(K_{\infty,v}) \otimes {\bf Q}_{p}/{\bf Z}_{p})$, and 
the surjectivity of the third map follows from 
Greenberg Lemma 4.6 in \cite{Gr2}.
Let $E[p^{\infty}]^{0}$ be the kernel of 
$E[p^{\infty}]=E(\overline{{\bf Q}})[p^{\infty}] \longrightarrow 
E(\overline{{\bf F}}_{p})[p^{\infty}]$ and 
$E[p^{\infty}]_{et}=E[p^{\infty}]/E[p^{\infty}]^{0}$. 
For a prime $v$ of $K_{\infty}$ above $p$,  
we denote by $K_{\infty,v,nr}$ the maximal unramified extension 
of $K_{\infty,v}$, and put 
$\Gamma_{v}=\Gal(K_{\infty,v,nr}/K_{\infty,v})$.
We know the isomorphism 
$H^{2}_{v}(K_{\infty,v}) \stackrel{\simeq}{\longrightarrow}
H^{1}(K_{\infty,v,nr}, E[p^{\infty}]_{et})^{\Gamma_{v}}$
by Greenberg \cite{Gr1} \S 2. 
If $v$ is a prime of $K_{\infty}$ not above $p$, we know 
$H^{2}_{v}(K_{\infty,v})=H^{1}(K_{\infty,v}, E[p^{\infty}])$. 
Therefore, we get an isomorphism 
$$(\bigoplus_{v \mid \frac{m'}{m}} H^{2}_{v}(K_{\infty,v}))^{G}
=\bigoplus_{u \mid \frac{m'}{m}} H^{2}_{u}({\bf Q}_{\infty,v})$$
where $v$ (resp. $u$) runs over all primes of $K_{\infty}$ 
(resp. ${\bf Q}_{\infty}$) 
above $m'/m=pN_{E}$. 
Thus, 
$\Sel(O_{K_{\infty}}[1/m'], E[p^{\infty}])^{G}
\longrightarrow 
\bigoplus_{v \mid \frac{m'}{m}} H^{2}_{v}(K_{\infty,v})^{G}$
is surjective. 
On the other hand, we have 
$H^{2}_{et}(\Spec O_{K_{\infty}}[1/m'],E[p^{\infty}])=0$ 
(see \cite{Gr1} Propositions 3, 4). 
This implies that 
$$H^{1}(G, H^{1}_{et}(\Spec O_{K_{\infty}}[1/m'],E[p^{\infty}]))
=H^{1}(G, \Sel(O_{K_{\infty}}[1/m'], E[p^{\infty}]))=0.$$ 
Taking the cohomology of the exact sequence (\ref{ES2}), 
we get 
\begin{equation} \label{Eq1}
H^{1}(G, \Sel(O_{K_{\infty}}[1/m], E[p^{\infty}]))=0.
\end{equation}
  
Therefore, $\Sel(O_{K_{\infty}}[1/m], E[p^{\infty}])$ 
is cohomologically trivial as a $G$-module by Serre 
\cite{CL} Chap. IX Th\'{e}or\`{e}me 8. 
This implies that 
$\Sel(O_{K_{\infty}}[1/m], E[p^{\infty}])^{\vee}$ is also 
cohomologically trivial. 
Since $\Sel(O_{K_{\infty}}[1/m], E[p^{\infty}])^{\vee}$
has no nontrivial finite submodule, 
the projective dimension of 
$\Sel(O_{K_{\infty}}[1/m], E[p^{\infty}])^{\vee}$
as a $\Lambda_{K_{\infty}}$-module
is $\leq 1$ by Popescu \cite{Pope} Proposition 2.3. 

\begin{Theorem} \label{T1}
Let $K/{\bf Q}$ be a finite abelian 
$p$-extension in which all bad primes of $E$ are unramified.
We take a finite set $S$ of good reduction primes which contains 
all ramifying primes in $K/{\bf Q}$ except $p$. 
Let $m$ be the product of primes in $S$.
Then \\
{\rm (1)} $\xi_{K_{\infty},S}$ is in the initial Fitting ideal 
$\Fitt_{0, \Lambda_{K_{\infty}}}
(\Sel(O_{K_{\infty}}[1/m], E[p^{\infty}])^{\vee})$. \\
{\rm (2)} We have 
$$\Fitt_{0, \Lambda_{K_{\infty}}}
(\Sel(O_{K_{\infty}}[1/m], E[p^{\infty}])^{\vee}) 
= \xi_{K_{\infty},S}\Lambda_{K_{\infty}}$$
if and only if the main conjecture (\ref{MC}) 
for $(E, {\bf Q}_{\infty}/{\bf Q})$ holds. 
\end{Theorem}

\noindent Proof. 
As we explained in the proof of Theorem \ref{GreenbergAMS}, 
we may assume that 
$K \cap {\bf Q}_{\infty} = {\bf Q}$. 
We recall that $\Sel(O_{K_{\infty}}[1/m], E[p^{\infty}])^{\vee}$ 
is a free ${\bf Z}_{p}$-module of finite rank 
under our assumptions. 

\vspace{5mm}

\noindent (1) 
Let $\psi:\Gal(K/{\bf Q}) \longrightarrow \overline{{\bf Q}}_{p}^{\times}$ 
be a character of $\Gal(K/{\bf Q})$, not necessarily faithful. 
We study the Fitting ideal of the $\psi$-quotient 
$(\Sel(O_{K_{\infty}}[1/m], E[p^{\infty}])^{\vee})_{\psi}= 
\Sel(O_{K_{\infty}}[1/m], E[p^{\infty}])^{\vee} \otimes_{{\bf Z}_{p}[\Gal(K/{\bf Q})]} O_{\psi}$. 
We denote by $F$ the subfield of $K$ corresponding to the 
kernel of $\psi$. 
We regard $\psi$ as a faithful character of $\Gal(F/{\bf Q})$.
Since $\Sel(O_{K_{\infty}}[1/m], E[p^{\infty}])^{\Gal(K/F)}
=\Sel(O_{F_{\infty}}[1/m], E[p^{\infty}])$, we have 
$$(\Sel(O_{K_{\infty}}[1/m], E[p^{\infty}])^{\vee})_{\psi}
=(\Sel(O_{F_{\infty}}[1/m], E[p^{\infty}])^{\vee})_{\psi}$$
where the right hand side is defined to be 
$\Sel(O_{F_{\infty}}[1/m], E[p^{\infty}])^{\vee} 
\otimes_{{\bf Z}_{p}[\Gal(F/{\bf Q})]} O_{\psi}$.

We put $\Lambda_{\psi}=(\Lambda_{F_{\infty}})_{\psi}$.
The group homomorphism $\psi$ induces the ring 
homomorphism $\Lambda_{F_{\infty}} \longrightarrow 
\Lambda_{\psi}$ which we also denote by $\psi$. 
The composition with $c_{K_{\infty}/F_{\infty}}:\Lambda_{K_{\infty}} 
\longrightarrow \Lambda_{F_{\infty}}$ and the above 
ring homomorphism $\psi$ is also 
denoted by $\psi: \Lambda_{K_{\infty}} \longrightarrow 
\Lambda_{\psi}$. 
Note that $F/{\bf Q}$ is a cyclic extension of degree a power of $p$. 
We denote by $F'$ the subfield of $F$ such that $[F:F']=p$. 
We put $N_{0}=N_{\Gal(F/F')}=\Sigma_{\sigma \in \Gal(F/F')} \sigma$.
If we put $[F:{\bf Q}]=p^{c}$ and take a generator $\gamma$ of $\Gal(F/{\bf Q})$, 
$N_{0}=\Sigma_{i=0}^{p-1} \gamma^{p^{c-1}i}$ is a cyclotomic polynomial 
and $O_{\psi}={\bf Z}_{p}[\mu_{p^{c}}] \simeq {\bf Z}_{p}[\Gal(F/{\bf Q})]/N_{0}$.
For any ${\bf Z}_{p}[\Gal(F/{\bf Q})]$-module $M$, 
we define $M^{\psi}=\Ker(N_{0}:M \longrightarrow M)$. 
Then 
the Pontrjagin dual of $M^{\psi}$ is 
$(M^{\psi})^{\vee}=(M^{\vee})/N_{0}=(M^{\vee})_{\psi}$. 
By the same method as the proof of (\ref{Eq1}), we have 
$H^{1}(\Gal(F/F'), \Sel(O_{F_{\infty}}[1/m], E[p^{\infty}]))=0$. 
Therefore, 
$\sigma-1: \Sel(O_{F_{\infty}}[1/m], E[p^{\infty}]) 
\longrightarrow 
\Sel(O_{F_{\infty}}[1/m], E[p^{\infty}])^{\psi}$
is surjective where $\sigma=\gamma^{p^{c-1}}$ is a generator of $\Gal(F/F')$. 
Therefore, taking the dual, we know that there is an injective 
homomorphism from 
$(\Sel(O_{F_{\infty}}[1/m], E[p^{\infty}])^{\vee})_{\psi}$ 
to 
$\Sel(O_{F_{\infty}}[1/m], E[p^{\infty}])^{\vee}$ which 
is a free ${\bf Z}_{p}$-module. 
Therefore, 
$(\Sel(O_{F_{\infty}}[1/m], E[p^{\infty}])^{\vee})_{\psi}$ 
contains no nontrivial finite $\Lambda_{\psi}$-submodule. 
This shows that 
$$\Fitt_{0, \Lambda_{\psi}}
((\Sel(O_{F_{\infty}}[1/m], E[p^{\infty}])^{\vee})_{\psi}) 
= \cc
((\Sel(O_{F_{\infty}}[1/m], E[p^{\infty}])^{\vee})_{\psi}).$$ 

Consider the $\psi$-quotient of the exact sequence (\ref{ES1});
$$
(\bigoplus_{v \mid m} T_{p}(E)_{\Gamma_{v}})_{\psi} 
\longrightarrow (\Sel(O_{F_{\infty}}[1/m], E[p^{\infty}])^{\vee})_{\psi}
\longrightarrow  (\Sel(O_{F_{\infty}}, E[p^{\infty}])^{\vee})_{\psi}
\longrightarrow  0
$$
where $v$ runs over all primes of $F_{\infty}$ above $m$. 
Since $\Ext_{{\bf Z}_{p}[\Gal(F/{\bf Q})]}^{1}
(O_{\psi}, \Sel(O_{F_{\infty}}, E[p^{\infty}])) 
={\hat H}^{0}(\Gal(F/{\bf Q}), \Sel(O_{F_{\infty}}, E[p^{\infty}]))$ 
is finite, the first map of the above exact sequence 
has finite kernel. 

Suppose that $\ell$ is a prime divisor of $m$. 
If $\ell$ is unramified in $F$, we have  
$$\Fitt_{0, \Lambda_{\psi}}
((\bigoplus_{v \mid \ell} T_{p}(E)_{\Gamma_{v}})_{\psi})
=P'_{\ell}(\sigma_{\ell})\Lambda_{\psi}$$
where $P'_{\ell}(x)=x^{2}-a_{\ell}x+\ell$. 
If $\ell$ is ramified in $F$, $\psi(\ell)=0$ and 
$(\bigoplus_{v \mid \ell} T_{p}(E)_{\Gamma_{v}})_{\psi}$ is finite. 
Therefore, we have 
$$\cc((\bigoplus_{v \mid m} T_{p}(E)_{\Gamma_{v}})_{\psi})
=(\prod_{\ell \in S \setminus S_{{\rm ram}}(F)}
P'_{\ell}(\sigma_{\ell}))\Lambda_{\psi}.
$$
Using the above exact sequence and Kato's theorem 
$\psi(\xi_{F_{\infty}}) \in 
\cc((\Sel(O_{F_{\infty}}, E[p^{\infty}])^{\vee})_{\psi})$, 
we have 
$$
\cc((\Sel(O_{F_{\infty}}[1/m], E[p^{\infty}])^{\vee})_{\psi})
\supset \psi(\xi_{F_{\infty}}) (\prod_{\ell \in S \setminus S_{{\rm ram}}(F)}
P'_{\ell}(\sigma_{\ell}))\Lambda_{\psi}.$$
Since $\xi_{F_{\infty}} (\prod_{\ell \in S \setminus S_{{\rm ram}}(F)}
P'_{\ell}(\sigma_{\ell}))=\xi_{F_{\infty},S}$ modulo unit 
and $c_{K_{\infty}/F_{\infty}}(\xi_{K_{\infty},S})
=\xi_{F_{\infty},S}$ by (\ref{Res4}), 
we obtain 
\begin{equation} \label{Inc1}
\psi(\xi_{K_{\infty},S}) \in 
\Fitt_{0, (\Lambda_{K_{\infty}})_{\psi}}
((\Sel(O_{K_{\infty}}[1/m], E[p^{\infty}])^{\vee})_{\psi})
\end{equation}
for any character $\psi$ of $\Gal(K/{\bf Q})$. 
Since the $\mu$-invariant of 
$\Sel(O_{K_{\infty}}[1/m], E[p^{\infty}])^{\vee}$ is zero 
as we explained above, (\ref{Inc1}) implies 
$$
\xi_{K_{\infty},S} \in 
\Fitt_{0, \Lambda_{K_{\infty}}}
(\Sel(O_{K_{\infty}}[1/m], E[p^{\infty}])^{\vee})
$$
(see Lemma 4.1 in \cite{Ku4}, for example).

\vspace{5mm}

\noindent (2)
We use the same notation $\psi$, $F$, etc. as above. 
At first, we assume (\ref{MC}). 
Then the algebraic $\lambda$-invariant of 
$\Sel(E/F_{\infty}, E[p^{\infty}])^{\vee}$ equals 
the analytic $\lambda$-invariant by 
Hachimori and Matsuno \cite{HM1}, \cite{HM2}, 
so the main conjecture 
$\cc((\Sel(O_{F_{\infty}}, E[p^{\infty}])^{\vee})_{\psi})
=\psi(\xi_{F_{\infty}}) \Lambda_{\psi}$ 
also holds. 
Therefore, we have 
\begin{eqnarray*} 
\cc((\Sel(O_{F_{\infty}}[1/m], E[p^{\infty}])^{\vee})_{\psi})
&=& \psi(\xi_{F_{\infty}}) (\prod_{\ell \in S \setminus S_{{\rm ram}}(F)}
P'_{\ell}(\sigma_{\ell}))\Lambda_{\psi}\\
&=&\psi(\xi_{F_{\infty},S})\Lambda_{\psi}
=\psi(\xi_{K_{\infty},S})\Lambda_{\psi}
\end{eqnarray*} 
and 
$$
\Fitt_{0, (\Lambda_{K_{\infty}})_{\psi}}
((\Sel(O_{K_{\infty}}[1/m], E[p^{\infty}])^{\vee})_{\psi})
=\psi(\xi_{K_{\infty},S}\Lambda_{K_{\infty}})\Lambda_{\psi}. 
$$
It follows from  \cite{Ku4} Corollary 4.2 that 
$$\Fitt_{0, \Lambda_{K_{\infty}}}
(\Sel(O_{K_{\infty}}[1/m], E[p^{\infty}])^{\vee})=
\xi_{K_{\infty},S}\Lambda_{K_{\infty}}.
$$

On the other hand, if we assume the above equality, 
taking the $\Gal(K/{\bf Q})$-invariant part of 
$\Sel(O_{K_{\infty}}[1/m], E[p^{\infty}])$, 
we get 
$$\Fitt_{0, \Lambda_{{\bf Q}_{\infty}}}
(\Sel(O_{{\bf Q}_{\infty}}[1/m], E[p^{\infty}])^{\vee})=
\xi_{{\bf Q}_{\infty},S}\Lambda_{{\bf Q}_{\infty},S},
$$
which 
implies (\ref{MC}).

\subsection{An analogue of Stickelberger's theorem} \label{subsec23}

Let $K/{\bf Q}$ be a finite abelian $p$-extension. 
When the conductor of $K$ is $m$, we define 
$\vartheta_{K} \in R_{K}={\bf Z}_{p}[\Gal(K/{\bf Q})]$ 
to be the image of 
$\vartheta_{{\bf Q}(\mu_{mp^{\infty}})}
\in \Lambda_{{\bf Q}(\mu_{mp^{\infty}})}$. 
Therefore, if $m$ is prime to $p$, 
$\vartheta_{K}$ is 
the image of 
$\vartheta_{{\bf Q}(\mu_{m})}=(1-\frac{\sigma_{p}}{\alpha})
(1-\frac{\sigma_{p}^{-1}}{\alpha})\tilde{\theta}_{{\bf Q}(\mu_{m})}$
by (\ref{Res3}). 
If $m=m'p^{n}$ for some $m'$ which is prime to $p$ and for some $n \geq 2$, 
$\vartheta_{K}$ is 
the image of 
$\vartheta_{{\bf Q}(\mu_{m'p^{n}})}=
\alpha^{-n}(\tilde{\theta}_{{\bf Q}(\mu_{m'p^{n}})}-
\alpha^{-1}\nu_{m'p^{n},m'p^{n-1}}
(\tilde{\theta}_{{\bf Q}(\mu_{m'p^{n-1}})}))$. 

For any positive integer $n$,  
we denote by ${\bf Q}(n)$ the maximal $p$-subextension 
of ${\bf Q}$ in ${\bf Q}(\mu_{n})$. 

\begin{Theorem} \label{AnnTh1}
For any finite abelian $p$-extension $K$ in which 
all bad primes of $E$ are unramified,
$\vartheta_{K}$ annihilates 
$\Sel(O_{K}, E[p^{\infty}])^{\vee}$, namely 
we have 
$$\vartheta_{K} \Sel(O_{K}, E[p^{\infty}])^{\vee} = 0.$$ 
\end{Theorem}

\noindent Proof. 
We may assume $K={\bf Q}(mp^{n})$ for some 
squarefree product $m$ of primes in ${\cal P}_{good}$ 
and for some $n \in {\bf Z}_{\geq 0}$. 
By Theorem \ref{T1} (1), taking $S$ to be the set of 
all prime divisors of $m$, 
we have $\xi_{K_{\infty}} \in 
\Fitt_{0, \Lambda_{K_{\infty}}}
(\Sel(O_{K_{\infty}}[1/m], E[p^{\infty}])^{\vee})$,  
which implies 
$\xi_{K_{\infty}} 
\Sel(O_{K_{\infty}}, E[p^{\infty}])^{\vee}=0$.
Let $\xi_{K} \in R_{K}={\bf Z}_{p}[\Gal(K/{\bf Q})]$ 
be the image of $\xi_{K_{\infty}}$. 
Since the natural map 
$\Sel(O_{K}, E[p^{\infty}]) \longrightarrow 
\Sel(O_{K_{\infty}}, E[p^{\infty}])$ is injective, 
we have 
$\xi_{K} \Sel(O_{K}, E[p^{\infty}])^{\vee}=0$.

By the definitions of $\xi_{{\bf Q}(\mu_{mp^{\infty}})}$, 
$\xi_{{\bf Q}(mp^{n})}$, 
$\vartheta_{{\bf Q}(mp^{n})}$,
we can write 
\begin{equation} \label{eq23}
\xi_{K}=
\xi_{{\bf Q}(mp^{n})}=\vartheta_{{\bf Q}(mp^{n})}+
\sum_{d \mid m, d \neq m}
\lambda_{d} \nu_{m,d}(\vartheta_{{\bf Q}(dp^{n})})
\end{equation} 
for some 
$\lambda_{d} \in R_{{\bf Q}(mp^{n})}$  
where $\nu_{m,d}: R_{{\bf Q}(dp^{n})} \longrightarrow 
R_{{\bf Q}(mp^{n})}$ 
is the norm map defined similarly as in \S \ref{subsec21}. 
We will prove this theorem by induction on $m$. 
Since $d<m$, we have 
$\vartheta_{{\bf Q}(dp^{n})} \in 
\Ann_{R_{{\bf Q}(dp^{n})}}
(\Sel(O_{{\bf Q}(dp^{n})}, E[p^{\infty}])^{\vee})$
by the hypothesis of the induction. 
This implies that 
$\nu_{m, d}(\vartheta_{{\bf Q}(dp^{n})})$ annihilates 
$\Sel(O_{{\bf Q}(mp^{n})}, E[p^{\infty}])^{\vee}$. 
Since $\xi_{K}$ is in 
$\Ann_{R_{K}}(\Sel(O_{K}, E[p^{\infty}])^{\vee})$, 
the above equation implies that $\vartheta_{K}$ is in 
$\Ann_{R_{K}}(\Sel(O_{K}, E[p^{\infty}])^{\vee})$.

\begin{Remark} \label{FittingIdeal00}
\begin{rm}
Let $K$, $S$, $m$ be as in Theorem \ref{T1}. 
Under our assumptions, the control theorem works completely; 
$$\Sel(O_{K}[1/m], E[p^{\infty}]) \stackrel{\simeq}{\longrightarrow}
\Sel(O_{K_{\infty}}[1/m], E[p^{\infty}])^{\Gal(K_{\infty}/K)}.$$
Therefore, 
Theorem \ref{T1} (1) implies that $\Fitt_{0, R_{K}}
(\Sel(O_{K}[1/m], E[p^{\infty}])^{\vee})$ is principal and 
\begin{equation} \label{FittingIdeal02}
\xi_{K,S} \in 
\Fitt_{0, R_{K}}
(\Sel(O_{K}[1/m], E[p^{\infty}])^{\vee})
\end{equation}
where $\xi_{K,S}$ is the image of $\xi_{K_{\infty},S}$ in $R_{K}$. 

Theorem \ref{T1} (2) implies that if we assume the main conjecture 
(\ref{MC}),  we have 
\begin{equation} \label{FittingIdeal0}
\Fitt_{0, R_{K}}
(\Sel(O_{K}[1/m], E[p^{\infty}])^{\vee})=\xi_{K,S}R_{K} \ .
\end{equation}
\end{rm}
\end{Remark}

\subsection{Higher Fitting ideals} \label{subsec24}

For a commutative ring $R$ and a finitely presented $R$-module $M$ 
with $n$ generators, let $A$ be an $n \times m$ relation matrix 
of $M$. 
For an integer $i \geq 0$, $\Fitt_{i,R}(M)$ is defined to be the 
ideal of $R$ generated by all $(n-i) \times (n-i)$ minors of $A$ 
(cf. \cite{N}; this ideal $\Fitt_{i,R}(M)$ does not depend on the choice 
of a relation matrix $A$).

Suppose that 
$K/{\bf Q}$ is a finite extension such that $K$ is in the 
cyclotomic ${\bf Z}_{p}$-extension ${\bf Q}_{\infty}$ of ${\bf Q}$, 
and that $m$ is a squarefree product of primes in 
${\cal P}^{(N)}$. 
We define $K(m)$ by $K(m)={\bf Q}(m)K$. 

We put 
${\cal G}_{\ell}=\Gal({\bf Q}(\ell)/{\bf Q})$ and
${\cal G}_{m}= \Gal({\bf Q}(m)/{\bf Q})
=\Pi_{\ell \mid m} {\cal G}_{\ell}$. 
We have $\Gal(K(m)/K)={\cal G}_{m}$. 
We put $n_{\ell}=\ord_{p}(\ell-1)$.
Suppose that $m=\ell_{1} \cdot...\cdot \ell_{r}$.
We take a generator $\tau_{\ell_{i}}$ of ${\cal G}_{\ell_{i}}$  
and put $S_{i}=\tau_{\ell_{i}}-1 \in R_{K(m)}$.
We write $n_{i}$ for $n_{\ell_{i}}$. 
We identify $R_{K(m)}$ with 
$$R_{K}[{\cal G}_{m}]=
R_{K}[S_{1},...,S_{r}]/
((1+S_{1})^{p^{n_{1}}}-1,...,
(1+S_{r})^{p^{n_{r}}}-1).$$
We consider $\vartheta_{K(m)} \in R_{K(m)}$ and write 
$$\vartheta_{K(m)}=\sum_{i_{1},...,i_{r} \geq 0}
a_{i_{1},...,i_{r}}^{(m)} S_{1}^{i_{i}}\cdot...\cdot S_{r}^{i_{r}}
$$
where $a_{i_{1},...,i_{r}}^{(m)} \in R_{K}$. 
Put $n_{0}=\min\{n_{1},...,n_{r}\}$.
For $s \in {\bf Z}_{>0}$, 
we define $c_{s}$ to be the maximal positive integer $c$  
such that 
$$T^{-1}((1+T)^{p^{n_{0}}}-1) 
\in p^{c}{\bf Z}_{p}[T]+ T^{s+1}{\bf Z}_{p}[T].$$
For example, $c_{1}=n_{0}$, ..., $c_{p-2}=n_{0}$,
$c_{p-1}=n_{0}-1$, ..., $c_{p^{2}-1}=n_{0}-2$.
If $i_{1}$,...,$i_{r} \leq s$, 
$a_{i_{1},...,i_{r}}^{(m)}$ mod $p^{c_{s}}$ 
is well-defined (it does not depend on the choice of 
$a_{i_{1},...,i_{r}}^{(m)}$). 

\begin{Theorem} \label{HigherFittingIdeal}
Let $K$ be an intermediate field of the cyclotomic 
${\bf Z}_{p}$-extension ${\bf Q}_{\infty}/{\bf Q}$ with 
$[K:{\bf Q}]<\infty$. 
Let $c_{s}$ be the integer defined above for $s \in {\bf Z}_{>0}$ 
and $m$.  
Assume that $i_{1}$,...,$i_{r} \leq s$ and 
$i_{1}+...+i_{r} \leq i$. 
Then we have 
$$a_{i_{1},...,i_{r}}^{(m)} \in 
\Fitt_{i, R_{K}/p^{c_{s}}}(\Sel(E/K, E[p^{c_{s}}])^{\vee}).$$
\end{Theorem}

For $m=\ell_{1} \cdot...\cdot \ell_{r}$, we denote $(-1)^{r}$ times 
the coefficient of $S_{1} \cdot ... \cdot S_{r}$ in 
$\vartheta_{K(m)}$ by $\delta_{m}$.
If $\ell_{i}$ splits completely in $K$ for all 
$i=1$,...,$r$, we can write 
\begin{equation} \label{Dcoeff0}
\vartheta_{K(m)}  \equiv 
\delta_{m}
\prod_{i=1}^{r}(1-\tau_{\ell_{i}}) 
=(-1)^{r}\delta_{m}S_{1}\cdot...\cdot S_{r}
\pmod{p^{N}, 
S_{1}^{2},...,S_{r}^{2}}
\end{equation}
(see \cite{Ku6} \S 6.3).
Taking $s=1$ and $i=r$ in Theorem \ref{HigherFittingIdeal}, we get  

\begin{Corollary} \label{HigherFittingIdeal2}
Let $K/{\bf Q}$ be a finite extension such that $K \subset {\bf Q}_{\infty}$. 
We have 
$$\delta_{m} \in \Fitt_{r, R_{K}/p^{N}}
(\Sel(E/K, E[p^{N}])^{\vee})$$
where $m=\ell_{1} \cdot...\cdot \ell_{r}$.
\end{Corollary}

\noindent Proof of Theorem \ref{HigherFittingIdeal}. 
We may assume $K={\bf Q}(p^{n})$ for some $n \geq 0$, 
so $K(m)={\bf Q}(mp^{n})$. 
First of all, we consider the image 
$\xi_{K(m)} \in R_{K(m)}$ of $\xi_{K(m)_{\infty}}$. 
Since $\Sel(E/K(m), E[p^{\infty}]) \longrightarrow 
\Sel(E/K(m)_{\infty}, E[p^{\infty}])$ is injective, 
$\xi_{K(m)}$ is in 
$\Fitt_{0, R_{K(m)}}(\Sel(E/K(m), E[p^{\infty}])^{\vee})$ 
by Theorem \ref{T1} (1). 
We write 
$$\xi_{K(m)}=\sum_{i_{1},...,i_{r} \geq 0}
\alpha_{i_{1},...,i_{r}}^{(m)} S_{1}^{i_{i}}\cdot...\cdot S_{r}^{i_{r}}
$$
where $\alpha_{i_{1},...,i_{r}}^{(m)} \in R_{K}$. 
Assume that $i_{1}$,...,$i_{r} \leq s$ and 
$i_{1}+...+i_{r} \leq i$. 
Then by Lemma 3.1.1 in \cite{Ku6} we have 
$$\alpha_{i_{1},...,i_{r}}^{(m)} \in 
\Fitt_{i, R_{K}/p^{c_{s}}}(\Sel(E/K, E[p^{c_{s}}])^{\vee}).$$

On the other hand, since $K(m)={\bf Q}(mp^{n})$ for some 
$n \geq 0$, we have 
$$\xi_{K(m)}=\vartheta_{K(m)}+
\sum_{d \mid m, d \neq m}
\lambda_{d} \nu_{m,d}(\vartheta_{{\bf Q}(dp^{n})})$$
for some $\lambda_{d} \in R_{K(m)}$  
by (\ref{eq23}). 
This implies that the images of $\xi_{K(m)}$ and $\vartheta_{K(m)}$ 
under the canonical homomorphism 
$$R_{K(m)}=R_{K}[S_{1},...,S_{r}]/I \longrightarrow 
R_{K}[[S_{1},...,S_{r}]]/J$$
coincide 
where $I=((1+S_{1})^{p^{n_{1}}}-1,...,(1+S_{1})^{p^{n_{r}}}-1)$  
and $J=(S_{1}^{-1}(1+S_{1})^{p^{n_{1}}}-1,...,
S_{r}^{-1}(1+S_{1})^{p^{n_{r}}}-1, S_{1}^{s+1},...,S_{r}^{s+1})$. 
Therefore, 
$\alpha_{i_{1},...,i_{r}}^{(m)} \equiv
a_{i_{1},...,i_{r}}^{(m)}$ mod $p^{c_{s}}$ for 
$i_{1}$,...,$i_{r} \leq s$. 
It follows that  
$a_{i_{1},...,i_{r}}^{(m)} \in 
\Fitt_{i, R_{K}/p^{c_{s}}}(\Sel(E/K, E[p^{c_{s}}])^{\vee})$.
This completes the proof of Theorem \ref{HigherFittingIdeal}.

\section{Review of Kolyvagin systems of Gauss sum type for 
elliptic curves} \label{section3}

In this section, we recall the results in \cite{Ku6} on 
Euler systems and Kolyvagin systems of Gauss sum type 
in the case of elliptic curves. 
From this section we assume all the assumptions (i), (ii), (iii), (iv) 
in \S \ref{subsec11}. 

\subsection{Some definitions} \label{subsec31}

Recall that in \S \ref{section2} we defined ${\cal P}_{good}$ by 
${\cal P}_{good}=\{\ell \mid $ $\ell$ is a 
good reduction prime for $E$ $\} \setminus \{p\}$, and ${\cal P}^{(N)}$ 
by
$${\cal P}^{(N)}=\{\ell \in {\cal P}_{good}
\mid \ell \equiv 1 \mbox{ (mod $p^{N}$)}\}$$
for a positive integer $N>0$. 
If $\ell$ is in ${\cal P}_{good}$, the absolute Galois group 
$G_{{\bf F}_{\ell}}$ acts on the group $E[p^{N}]$ of $p^{N}$-torsion 
points, so 
we consider $H^{i}({\bf F}_{\ell}, E[p^{N}])$. 
We define 
\begin{eqnarray*}
{\cal P}_{0}^{(N)}&=&\{\ell \in {\cal P}^{(N)} \mid 
H^{0}({\bf F}_{\ell}, E[p^{N}]) \ 
\mbox{contains an element of order $p^{N}$}\},\\
({\cal P}'_{0})^{(N)}&=&\{\ell \in {\cal P}^{(N)} \mid 
H^{0}({\bf F}_{\ell}, E[p^{N}])=E[p^{N}]\}, \ \mbox{and}\\
{\cal P}_{1}^{(N)}&=&\{\ell \in {\cal P}^{(N)} \mid 
H^{0}({\bf F}_{\ell}, E[p^{N}]) 
\simeq {\bf Z}/p^{N} \}.
\end{eqnarray*}
So ${\cal P}_{0}^{(N)} \supset ({\cal P}'_{0})^{(N)}$, 
${\cal P}_{0}^{(N)} \supset{\cal P}_{1}^{(N)}$, and  
$({\cal P}'_{0})^{(N)} \cap {\cal P}_{1}^{(N)}=\emptyset$. 
Suppose that $\ell$ is in ${\cal P}_{1}^{(N)}$. 
Then, since $\ell \equiv 1$ (mod $p^{N})$, we have an exact sequence 
$0 \longrightarrow {\bf Z}/p^{N} \longrightarrow E[p^{N}] 
\longrightarrow {\bf Z}/p^{N} \longrightarrow 0$ of 
$G_{{\bf F}_{\ell}}$-modules where 
$G_{{\bf F}_{\ell}}$ acts on ${\bf Z}/p^{N}$ trivially. 
So the action of the Frobenius $\F_{\ell}$ at $\ell$ on $E[p^{N}]$ 
can be written as  
$
\left( \begin{array}{cc}
1 & 1\\
0 & 1
\end{array} \right)
$
for a suitable basis of $E[p^{N}]$.
Therefore, $H^{1}({\bf F}_{\ell}, E[p^{N}])$ is also isomorphic to 
${\bf Z}/p^{N}$ for $\ell \in {\cal P}_{1}^{(N)}$. 

Let $t \in E[p^{N}]$ be an element of order $p^{N}$. 
We define 
\begin{eqnarray*}
{\cal P}_{0,t}^{(N)}&=&\{\ell \in {\cal P}^{(N)} \mid 
t \in H^{0}({\bf F}_{\ell}, E[p^{N}])\}, \\ 
{\cal P}_{1,t}^{(N)}&=&\{\ell \in {\cal P}^{(N)} \mid 
H^{0}({\bf F}_{\ell}, E[p^{N}]) = ({\bf Z}/p^{N})t \}.
\end{eqnarray*}
So, 
${\cal P}_{0}^{(N)}=\bigcup_{t} {\cal P}_{0,t}^{(N)}$ and 
${\cal P}_{1}^{(N)}=\bigcup_{t} {\cal P}_{1,t}^{(N)}$
where $t$ runs over all elements of order $p^{N}$. 
Since we assumed that the Galois action on the Tate 
module is surjective, both   
$({\cal P}'_{0})^{(N)}$ and ${\cal P}_{1,t}^{(N)}$ 
are infinite by Chebotarev density theorem 
(\cite{Ku6} \S 4.3). 
                       
We define ${\cal K}_{(p)}$ to be the set of number fields $K$ such that 
$K/{\bf Q}$ is a finite abelian $p$-extension in which 
all bad primes of $E$ are unramified. 
Suppose that $K$ is in ${\cal K}_{(p)}$. 
We define 
\begin{eqnarray*}
({\cal P}'_{0})^{(N)}(K)&=&\{\ell \in ({\cal P}'_{0})^{(N)} \mid 
\mbox{$\ell$ splits completely in $K$}\},\\
{\cal P}_{1}^{(N)}(K)&=&\{\ell \in {\cal P}_{1}^{(N)} \mid 
\mbox{$\ell$ splits completely in $K$}\}.
\end{eqnarray*}
Again by Chebotarev density theorem, both 
$({\cal P}'_{0})^{(N)}(K)$ and ${\cal P}_{1}^{(N)}(K)$ 
are infinite (see \cite{Ku6} \S 4.3). 

\vspace{5mm}

Suppose $\ell \in {\cal P}_{good}$. 
For a prime $v$ above $\ell$, we know 
$H^{1}(K_{v}, E[p^{N}])/(E(K_{v}) \otimes {\bf Z}/p^{N}) 
=H^{0}(\kappa(v), E[p^{N}](-1))$ where 
$\kappa(v)$ is the residue field of $v$. 
We put 
\begin{equation} \label{H2ell}
{\cal H}_{\ell}^{2}(K)=\bigoplus_{v \mid \ell} 
H^{0}(\kappa(v), E[p^{N}](-1)).
\end{equation}
If $\ell$ is in $({\cal P}'_{0})^{(N)}(K)$ 
(resp. ${\cal P}_{1}^{(N)}(K)$), 
${\cal H}_{\ell}^{2}(K)$ is a free $R_{K}/p^{N}$-module 
of rank $2$ (resp. rank $1$) 
where $R_{K}={\bf Z}_{p}[\Gal(K/{\bf Q})]$ as before.

\vspace{5mm}

From now on, for a prime $\ell \in {\cal P}_{0}^{(N)}$, we fix a prime 
$\ell_{\overline{{\bf Q}}}$ of an algebraic closure 
$\overline{{\bf Q}}$ above $\ell$. 
For any algebraic number field $F$, 
we denote the prime of $F$ below 
$\ell_{\overline{{\bf Q}}}$ by $\ell_{F}$, 
so when we consider finite extensions $F_{1}/k$, $F_{2}/k$ such that 
$F_1 \subset F_2$, 
the primes $\ell_{F_{2}}$, $\ell_{F_{1}}$ 
satisfy $\ell_{F_{2}} \vert \ell_{F_{1}}$.

We take a primitive $p^{n}$-th root of unity $\zeta_{p^{n}}$ 
such that $(\zeta_{p^{n}})_{n \geq 1} \in {\bf Z}_{p}(1)=
{\lim\limits_{\longleftarrow}} \mu_{p^{n}}$, and fix it. 

\vspace{5mm}

In the following, for each $\ell$ in ${\cal P}_{0}^{(N)}(K)$, 
we take $t_{\ell} 
\in H^{0}({\bf F}_{\ell}, E[p^{N}])$ and fix it. 
We define 
\begin{equation}
t_{\ell,K}=(t_{\ell} \otimes \zeta_{p^{N}}^{\otimes (-1)},0,...,0) 
\in {\cal H}_{\ell}^{2}(K)
\end{equation}
where the right hand side is the element whose $\ell_{K}$-component is 
$t_{\ell} \otimes \zeta_{p^{N}}^{\otimes (-1)}$ and 
other components are zero. 

\vspace{5mm}

Suppose that $K$ is in ${\cal K}_{(p)}$. 
Let $K_{\infty}/K$ be the cyclotomic ${\bf Z}_{p}$-extension, 
and $K_{n}$ be the $n$-th layer. 
Since $\Sel(O_{K_{\infty}}, E[p^{\infty}])^{\vee}$
is a finitely generated ${\bf Z}_{p}$-module, 
the corestriction map 
$\Sel(O_{K_{m}}, E[p^{N}]) \longrightarrow 
\Sel(O_{K}, E[p^{N}])$ is the zero map 
if $m$ is sufficiently large.  
We take the minimal $m>0$ satisfying this property, and
put $K_{[1]}=K_{m}$. 
We define inductively $K_{[n]}$ by 
$K_{[n]}=(K_{[n-1]})_{[1]}$ 
where we applied the above definition to $K_{[n-1]}$ 
instead of $K$. 

We can compute how large $K_{[n]}$ is. 
Let $\lambda$ be the $\lambda$-invariant of 
$\Sel(O_{K_{\infty}}, E[p^{\infty}])^{\vee}$.
We take $a \in {\bf Z}_{\geq 0}$ such that 
$p^{a+1}-p^{a} \geq \lambda$. 
Suppose that $K=K'_{m}$ ($m$-th layer of $K'_{\infty}/K'$) 
for some $K'$ such that $p$ is unramified in $K'$.
The corestriction map 
$\Sel(O_{K'_{a+1}}, E[p]) \longrightarrow \Sel(O_{K'_{a}}, E[p])$ 
is the zero map. 
Therefore, 
$\Sel(O_{K'_{a+N}}, E[p^{N}]) \longrightarrow \Sel(O_{K'_{a}}, E[p^{N}])$ 
is the zero map. 
Put $a'=\max(a-m, 0)$.
Then $\Sel(O_{K_{a'+N}}, E[p^{N}]) \longrightarrow 
\Sel(O_{K_{a'}}, E[p^{N}])$ 
is the zero map.
Therefore, we have $K_{[1]} \subset K_{a'+N}$. 
Also we know $K_{[n]} \subset K_{a'+nN}$.

Let $n_{\lambda}$, $d_{n}$ be the numbers defined just before 
(\ref{PNn}) in \S \ref{subsec12}. 
Then we can show that if $\ell \in {\cal P}_{1}^{(N)}$ 
satisfies $\ell \equiv 1$ (mod $p^{d_{n}}$), 
$\ell$ is in ${\cal P}_{1}^{(N)}({\bf Q}_{[n]})$ by 
the same method as above.

\subsection{Euler systems of Gauss sum type for 
elliptic curves} \label{subsec32}

We use the following lemma which is the global duality 
theorem (see Theorem 2.3.4 in Mazur and Rubin \cite{MR}).

\begin{Lemma} \label{FES0}
Suppose that 
$m$ is a product of primes in ${\cal P}_{good}$. 
We have an exact sequence 
$$
0 \longrightarrow \Sel(O_{K}, E[p^{N}])  
\longrightarrow \Sel(O_{K}[1/m], E[p^{N}]) \longrightarrow 
\bigoplus_{\ell \mid m} {\cal H}_{\ell}^{2}(K) 
\longrightarrow \Sel(O_{K}, E[p^{N}])^{\vee}.  
$$
\end{Lemma}

We remark that we can take $m$ such that the last map is 
surjective in our case (see Lemma \ref{Lemma621} below). 

\vspace{5mm}

Let $K$ be a number field in ${\cal K}_{(p)}$ and 
$\ell \in {\cal P}_{0}^{(N)}(K_{[1]})$. 
We apply the above lemma to $K_{[1]}$ and obtain an 
exact sequence
$$\Sel(O_{K_{[1]}}[1/\ell], E[p^{N}]) 
\stackrel{\partial_{\ell}}{\longrightarrow} 
{\cal H}_{\ell}^{2}(K_{[1]})
\stackrel{w_{\ell}}{\longrightarrow} 
\Sel(O_{K_{[1]}}, E[p^{N}])^{\vee}.$$
Consider $\vartheta_{K_{[1]}} t_{\ell,K_{[1]}} \in {\cal H}_{\ell}^{2}(K_{[1]})$.
By Theorem \ref{AnnTh1} we know 
$w_{\ell}(\vartheta_{K_{[1]}} t_{\ell,K_{[1]}})=
\vartheta_{K_{[1]}}w_{\ell}(t_{\ell,K_{[1]}})
=0$. 
Therefore, there is an element $g \in \Sel(O_{K_{[1]}}[1/\ell], E[p^{N}])$ 
such that $\partial_{\ell}(g)= \vartheta_{K_{[1]}} t_{\ell,K_{[1]}}$. 
We define 
\begin{equation} \label{ES31}
g_{\ell,t_{\ell}}^{(K)} = \Cor_{K_{[1]}/K}(g) 
\in \Sel(O_{K}[1/\ell], E[p^{N}]). 
\end{equation}
This element $g_{\ell,t_{\ell}}^{(K)}$ does not depend on the 
choice of $g \in \Sel(O_{K_{[1]}}[1/\ell], E[p^{N}])$ 
(\cite{Ku6} \S 5.4).
We write $g_{\ell}$ instead of $g_{\ell,t_{\ell}}^{(K)}$
when no confusion arises. 

\begin{Remark}
\begin{rm}
To define $g_{\ell}$, we used in \cite{Ku6} the $p$-adic $L$-function 
$\theta_{K_{\infty}}$ whose 
Euler factor at $\ell$ is 
$1-\frac{a_{\ell}}{\ell} \sigma_{\ell}^{-1}+\frac{1}{\ell}
\sigma_{\ell}^{-2}$. 
The element $\theta_{K_{\infty}}$ can be constructed from 
$\vartheta_{K_{\infty}}$ by the same method as when we constructed 
$\xi_{K_{\infty}}$ in \S \ref{subsec21}. 
In the above definition (\ref{ES31}),
we used $\vartheta_{K}$ (namely $\vartheta_{K_{\infty}}$) 
instead of $\theta_{K_{\infty}}$.
\end{rm}
\end{Remark}

\subsection{Kolyvagin derivatives of Gauss sum type} \label{subsec33}

Let $\ell$ be a prime in ${\cal P}_{good}$.
We define $\partial_{\ell}$ as a natural homomorphism 
$$\partial_{\ell}: H^{1}(K, E[p^{N}]) \longrightarrow 
{\cal H}_{\ell}^{2}(K)=
\bigoplus_{v \mid \ell} H^{0}(\kappa(v), E[p^{N}](-1))$$
where we used 
$H^{1}(K_{v}, E[p^{N}])/(E(K_{v}) \otimes {\bf Z}/p^{N}) 
=H^{0}(\kappa(v), E[p^{N}](-1))$.

Next, we assume $\ell \in {\cal P}_{1}^{(N)}(K)$. 
We denote by ${\bf Q}_{\ell}(\ell)$ the maximal $p$-subextension 
of ${\bf Q}_{\ell}$ inside ${\bf Q}_{\ell}(\mu_{\ell})$.
Put ${\cal G}_{\ell}=\Gal({\bf Q}_{\ell}(\ell)/{\bf Q}_{\ell})$. 
By Kummer theory, ${\cal G}_{\ell}$ is isomorphic to 
$\mu_{p^{n_{\ell}}}$ where $n_{\ell}=\ord_{p}(\ell-1)$. 
We denote by $\tau_{\ell}$ the corresponding element 
of ${\cal G}_{\ell}$ to $\zeta_{p^{n_{\ell}}}$ that is the  
primitive $p^{n_{\ell}}$-th root of unity we fixed.

We consider the natural homomorphism 
$H^{1}({\bf Q}_{\ell}, E[p^{N}]) \longrightarrow
H^{1}({\bf Q}_{\ell}(\ell), E[p^{N}])$ 
and denote the kernel by $H^{1}_{tr}({\bf Q}_{\ell}, E[p^{N}])$.
Let ${\bf Q}_{\ell, nr}$ be the maximal unramified extension 
of ${\bf Q}_{\ell}$. 
We identify 
$H^{1}({\bf F}_{\ell}, E[p^{N}])$ 
with  
$H^{1}(\Gal({\bf Q}_{\ell, nr}/{\bf Q}_{\ell}), E[p^{N}])$, 
and regard it as a subgroup of 
$H^{1}({\bf Q}_{\ell}, E[p^{N}])$. 
Then both $H^{1}({\bf F}_{\ell}, E[p^{N}])$ and 
$H^{1}_{tr}({\bf Q}_{\ell}, E[p^{N}])$ are isomorphic to 
${\bf Z}/p^{N}$, and 
we have decomposition 
$$H^{1}({\bf Q}_{\ell}, E[p^{N}])=H^{1}({\bf F}_{\ell}, E[p^{N}])
\oplus H^{1}_{tr}({\bf Q}_{\ell}, E[p^{N}])$$
as an abelian group. 
We also note that $H^{1}({\bf F}_{\ell}, E[p^{N}])$ coincides with 
the image of the Kummer map and is isomorphic to 
$E({\bf Q}_{\ell}) \otimes {\bf Z}/p^{N}$. 
We consider the homomorphism 
\begin{eqnarray} \label{phiprime}
\phi': H^{1}({\bf Q}_{\ell}, E[p^{N}]) 
\longrightarrow
H^{1}({\bf F}_{\ell}, E[p^{N}]) 
\end{eqnarray} 
which is obtained from the above decomposition. 

Note that $H^{1}({\bf F}_{\ell}, E[p^{N}])=
E[p^{N}]/(\F_{\ell}-1)$ where $\F_{\ell}$ is 
the Frobenius at $\ell$. 
Since $\ell$ is in ${\cal P}_{1}^{(N)}$, 
$\F_{\ell}^{-1}-1:E[p^{N}]/(\F_{\ell}-1)
\longrightarrow E[p^{N}]^{\F_{\ell}=1}=
H^{0}({\bf F}_{\ell}, E[p^{N}])$ is an isomorphism. 
We define $\phi'': H^{1}({\bf Q}_{\ell}, E[p^{N}]) 
\longrightarrow H^{0}({\bf F}_{\ell}, E[p^{N}])$ 
as the composition of $\phi'$ and 
$H^{1}({\bf F}_{\ell}, E[p^{N}]) 
 \stackrel{\F_{\ell}^{-1}-1}{\longrightarrow}
H^{0}({\bf F}_{\ell}, E[p^{N}])$.
We define 
$$\phi_{\ell}:H^{1}(K, E[p^{N}]) \longrightarrow 
{\cal H}_{\ell}^{2}(K)(1)$$
as the composition of  
the natural homomorphism
$H^{1}(K, E[p^{N}]) \longrightarrow 
\bigoplus_{v \mid \ell} H^{1}(K_{v}, E[p^{N}])$
and $\phi''$ for $K_{v}$.
Using the primitive $p^{N}$-th root of unity $\zeta_{p^{N}}$ 
we fixed, we regard $\phi_{\ell}$ as a homomorphism 
$$\phi_{\ell}:H^{1}(K, E[p^{N}]) \longrightarrow 
{\cal H}_{\ell}^{2}(K).$$

\vspace{5mm}

For a prime $\ell \in {\cal P}_{1}^{(N)}(K)$, we put 
${\cal G}_{\ell}=\Gal({\bf Q}(\ell)/{\bf Q})$. 
We identify ${\cal G}_{\ell}$ with 
$\Gal({\bf Q}_{\ell}(\ell)/{\bf Q}_{\ell})$. 
Recall that we defined $n_{\ell}$ by 
$p^{n_{\ell}}=[{\bf Q}(\ell):{\bf Q}]$, and 
we took a generator $\tau_{\ell}$ of ${\cal G}_{\ell}$ above. 
We define 
$$N_{\ell}=\sum_{i=0}^{p^{n_{\ell}}-1} \tau_{\ell}^{i} 
\in {\bf Z}[{\cal G}_{\ell}], \ 
D_{\ell} = \sum_{i=0}^{p^{n_{\ell}}-1} i \tau_{\ell}^{i} 
\in {\bf Z}[{\cal G}_{\ell}]$$
as usual.  

We define ${\cal N}_{1}^{(N)}(K)$ to be the set of 
squarefree products of primes in ${\cal P}_{1}^{(N)}(K)$. 
We suppose $1 \in  {\cal N}_{1}^{(N)}(K)$.
For $m \in {\cal N}_{1}^{(N)}(K)$,
we put ${\cal G}_{m}=\Gal({\bf Q}(m)/{\bf Q})$, 
$N_{m}=\Pi_{\ell \vert m}N_{\ell} \in {\bf Z}[{\cal G}_{m}]$, and 
$D_{m}=\Pi_{\ell \vert m}D_{\ell} \in {\bf Z}[{\cal G}_{m}]$. 
Assume that $\ell$ is in 
$({\cal P}'_{0})^{(N)}(K(m)_{[1]})$ and consider 
$g_{\ell,t_{\ell}}^{K(m)}\in 
\Sel(O_{K(m)}[1/\ell], E[p^{N}])$. 
We can check that 
$D_{m} g_{\ell,t_{\ell}}^{K(m)}$ is in 
$\Sel(O_{K(m)}[1/m \ell], E[p^{N}])^{{\cal G}_{m}}$.
Using the fact that $\Sel(O_{K}[1/m \ell], E[p^{N}]) 
\stackrel{\simeq}{\longrightarrow}
\Sel(O_{K(m)}[1/m \ell], E[p^{N}])^{{\cal G}_{m}}$ 
is bijective by Lemma \ref{Lemma331} below 
(cf. also \cite{Ku6} Lemma 6.3.1), 
we define   
\begin{equation} \label{KDer}
\kappa_{m, \ell}=\kappa_{m, \ell, t_{\ell}}^{(K)} 
\in \Sel(O_{K}[1/m \ell], E[p^{N}])
\end{equation}
to be the unique element whose 
image in $\Sel(O_{K(m)}[1/m \ell], E[p^{N}])$ 
is $D_{m} g_{\ell,t_{\ell}}^{(K(m))}$.

The following lemma will be also used in the next section. 

\begin{Lemma} \label{Lemma331}
Suppose that $K$, $L \in {\cal K}_{(p)}$ and $K \subset L$. 
For any $m \in {\bf Z}_{>0}$, the restriction map 
$\Sel(O_{K}[1/m], E[p^{N}]) 
\stackrel{\simeq}{\longrightarrow} 
\Sel(O_{L}[1/m], E[p^{N}])^{\Gal(L/K)}$ 
is bijective.  
\end{Lemma}

\noindent Proof. Let $N_{E}$ be the conductor of $E$, 
$m'=mpN_{E}$, and $m''$ the product of primes 
which divide $pN_{E}$ and which do not divide $m$. 
Put $G=\Gal(L/K)$. 
We have a commutative diagram of exact sequences
$$
\begin{array}{ccccccc}
0 & \longrightarrow  & 
\Sel(O_{K}[1/m], E[p^{N}]) & \longrightarrow &
\Sel(O_{K}[1/m'], E[p^{N}]) &
\longrightarrow &
\bigoplus_{v \mid m''} 
H^{2}_{K,v} \\
&&\mapdown{\alpha_{1}}&&\mapdown{\alpha_{2}}&&
\mapdown{\alpha_{3}} \\
0 & \longrightarrow  & 
\Sel(O_{L}[1/m], E[p^{N}])^{G} & \longrightarrow &
\Sel(O_{L}[1/m'], E[p^{N}])^{G} &
\longrightarrow &
(\bigoplus_{w \mid m''} 
H^{2}_{L,w})^{G}
\end{array}
$$
where $H^{2}_{K,v}=H^{1}(K_{v}, E[p^{N}])/(E(K_{v}) 
\otimes {\bf Z}/p^{N})$ and 
$H^{2}_{L,w}=H^{1}(L_{w}, E[p^{N}])/(E(L_{w}) 
\otimes {\bf Z}/p^{N})$.
Since $\Sel(O_{L}[1/m'], E[p^{N}])=H^{1}_{et}
(\Spec O_{L}[1/m'], E[p^{N}])$ and 
$H^{0}(L, E[p^{N}])=0$, 
$\alpha_{2}$ is bijective. 
Suppose that $v$ divides $m''$ and $w$ is above $v$. 
When $v$ divides $N_{E}$, 
since $v$ is unramified in $L$ and $p$ is prime to $\Tam(E)$, 
$H^{2}_{K,v} \longrightarrow H^{2}_{L,w}$
is injective (Greenberg \cite{Gr2} \S 3). 
When $v$ is above $p$, 
$H^{2}_{K,v} \longrightarrow H^{2}_{L,w}$
is injective because $a_{p} \not \equiv 1$ (mod $p$) 
(Greenberg \cite{Gr2} \S 3). 
Hence $\alpha_{3}$ is injective.  
Therefore, $\alpha_{1}$ is bijective. 

\vspace{5mm}

In \cite{Ku5}, if 
$m$ has a factorization 
$m=\ell_{1}\cdot ...\cdot \ell_{r}$ such that 
$\ell_{i+1} \in {\cal P}_{1}^{(N)}(K( \ell_{1}\cdot ...\cdot \ell_{i}))$ 
for all $i=1$,...,$r-1$, we called $m$ {\it well-ordered}. 
But the word ``well-ordered" might cause confusion, so  
we call $m$ {\it admissible} in this paper 
if $m$ satisfies the above condition. 
Note that we do not impose the condition 
$\ell_{1}< ...< \ell_{r}$ in the above definition, and that 
$m$ is admissible if there is one factorization as above.
We sometimes call the set of prime divisors of $m$ {\it admissible} 
if $m$ is admissible.  

Suppose that $m=\ell_{1}\cdot...\cdot \ell_{r}$. 
We define  
$\delta_{m} \in R_{K}/p^{N}$ by 
\begin{equation} \label{Dcoeff}
\vartheta_{K(m)}  \equiv 
\delta_{m}
\prod_{i=1}^{r}(1-\tau_{\ell_{i}}) \pmod{p^{N}, 
(\tau_{\ell_{1}}-1)^{2},...,(\tau_{\ell_{r}}-1)^{2} }
\end{equation}
 (see \cite{Ku6} \S 6.3). 

We simply write $\kappa_{m, \ell}$ for $\kappa_{m, \ell, t_{\ell}}^{(K)}$.  
We have the following Proposition (\cite{Ku6} Propositions 6.3.2, 
6.4.5 and Lemma 6.3.4). 

\begin{Proposition} \label{KKP1}  
Suppose that $m$ is in ${\cal N}_{1}^{(N)}(K)$, and 
$\ell \in ({\cal P}'_{0})^{(N)}(K(m)_{[1]})$. 
We take $n_{0}$ sufficiently large such that 
every prime of $K_{n_{0}}$ dividing $m$ is inert in 
$K_{\infty}/K_{n_{0}}$.
We further assume that $\ell \in 
({\cal P}'_{0})^{(N)}(K_{n_{0}+N})$. 
Then \\
{\rm (0)} $\kappa_{m, \ell} \in 
\Sel(O_{K}[1/m \ell], E[p^{N}])$.\\
{\rm (1)} $\partial_{r}(\kappa_{m,\ell})=
\phi_{r}(\kappa_{\frac{m}{r},\ell})$ for any prime divisor 
$r$ of $m$. \\
{\rm (2)} $\partial_{\ell}(\kappa_{m,\ell})=\delta_{m}t_{\ell,K}$. \\
{\rm (3)} Assume further that $m$ is admissible. 
Then $\phi_{r}(\kappa_{m, \ell})=0$ for any prime divisor 
$r$ of $m$. 
\end{Proposition}

\subsection{Construction of Kolyvagin systems 
of Gauss sum type} \label{subsec34}

In the previous subsection we constructed $\kappa_{m,\ell}$ for 
$m \in {\cal N}_{1}^{(N)}(K)$ and a prime 
$\ell \in ({\cal P}'_{0})^{(N)}(K)$ satisfying 
some properties.
In this subsection we construct 
$\kappa_{m,\ell}$ for $\ell \in {\cal P}_{1}^{(N)}(K)$ 
satisfying some properties (see Proposition \ref{KKP2}).
The property (4) in Proposition \ref{KKP2} is a beautiful property
of our Kolyvagin systems of Gauss sum type, which is unique for 
Kolyvagin systems of Gauss sum type.

For a squarefree product $m$ of primes, 
we define $\epsilon(m)$ to be the number of prime divisors of 
$m$, namely $\epsilon(m)=r$ if $m=\ell_{1}\cdot...\cdot \ell_{r}$. 

For any prime number $\ell$, we write 
${\cal H}_{\ell}^{2}(K)=\bigoplus_{v \mid \ell} 
H^{1}(K_{v}, E[p^{N}])/(E(K_{v}) \otimes {\bf Z}/p^{N})$, 
and consider the natural map 
$$w_{K}: \bigoplus_{\ell} {\cal H}_{\ell}^{2}(K) 
\longrightarrow \Sel(O_{K}, E[p^{N}])^{\vee}$$
which is obtained by taking the dual of 
$\Sel(O_{K}, E[p^{N}]) \longrightarrow \bigoplus_{v} E(K_{v}) 
\otimes {\bf Z}/p^{N}$.
We also consider the natural map 
$$\partial_{K}: H^{1}(K, E[p^{N}])
\longrightarrow 
\bigoplus_{\ell} {\cal H}_{\ell}^{2}(K).$$ 
 
We use the following lemma which was proved in \cite{Ku6} 
Proposition 4.4.3 and Lemma 6.2.1 (2).

\begin{Lemma} \label{Lemma621}
Suppose that $K \in {\cal K}_{(p)}$ and  
$r_{1}$,...,$r_{s}$ are $s$ distinct primes in 
${\cal P}_{1}^{(N)}(K)$. 
Assume that for each $i=1$,...,$s$, 
$\sigma_{i} \in {\cal H}_{r_{i}}^{2}(K)$ is given, 
and also $x \in \Sel(O_{K}, E[p^{N}])^{\vee}$ is given. 
Let $K'/K$ be an extension such that 
$K' \in {\cal K}_{(p)}$. 
Then there are infinitely many $\ell \in  {\cal P}_{0}^{(N)}(K)$ 
such that $w_{K}(t_{\ell,K})=x$.
We take such a prime $\ell$ and fix it. 
Then there are infinitely many $\ell' \in  ({\cal P}'_{0})^{(N)}(K')$ 
which satisfy the following properties: \\
$({\rm i})$ 
$w_{K}(t_{\ell',K})=w_{K}(t_{\ell,K})=x$. \\
$({\rm ii})$ There is an element $z \in 
\Sel(O_{K}[1/\ell \ell'], E[p^{N}])$ 
such that 
$\partial_{K}(z) = t_{\ell',K}-t_{\ell,K}$ and 
$\phi_{r_{i}}(z)=\sigma_{i}$ for each 
$i=1$,...,$s$. 
\end{Lemma}

Assume that $m\ell$ is in 
${\cal N}_{1}^{(N)}(K_{[\epsilon(m\ell)]})$. 
By Lemma \ref{Lemma621} we can take 
$\ell' \in  ({\cal P}'_{0})^{(N)}$ satisfying the 
following properties: \\
(i) $\ell' \in  ({\cal P}'_{0})^{(N)}(K_{[\epsilon(m\ell)]}(m)_{[1]}K_{n_{0}+N})$ 
where $n_{0}$ is as in Proposition \ref{KKP1}. \\
(ii) $w_{K_{[\epsilon(m\ell)]}}(t_{\ell', K_{[\epsilon(m\ell)]}})
=w_{K_{[\epsilon(m\ell)]}}(t_{\ell, K_{[\epsilon(m\ell)]}})$. \\
(iii) Let $\phi_{r}^{(K_{[\epsilon(m\ell)]})}: 
H^{1}(K_{[\epsilon(m\ell)]}, E[p^{N}]) \longrightarrow 
{\cal H}_{r}^{2}(K_{[\epsilon(m\ell)]})$ be the map 
$\phi_{r}$ for $K_{[\epsilon(m\ell)]}$. 
There is an element $b'$ in 
$\Sel(O_{K_{[\epsilon(m\ell)]}}[1/\ell \ell'], E[p^{N}])$ 
such that 
$$
\partial_{K_{[\epsilon(m \ell)]}}(b')=
t_{\ell', K_{[\epsilon(m \ell)]}} - t_{\ell, K_{[\epsilon(m \ell)]}}$$
and  $\phi_{r}^{K_{[\epsilon(m \ell)]}}(b')=0$ for all $r$ dividing $m$. 

We have already defined $\kappa_{m,\ell'}$ in the previous subsection. 
We put $b=\Cor_{K_{[\epsilon(m \ell)]}/K}(b')$ and define 
\begin{equation}
\kappa_{m,\ell}=\kappa_{m,\ell'}- \delta_{m} b.
\end{equation}
Then this element does not depend on the choice of $\ell'$ and $b'$ 
(see \cite{Ku6} \S 6.4). 
In \cite{Ku6}, we took $b'$ which does not necessarily satisfy 
$\phi_{r}^{K_{[\epsilon(m \ell)]}}(b')=0$ in the definition of 
$\kappa_{m,\ell}$. 
But we adopted the above definition here because it is simpler and 
there is no loss of generality.  
 
The next proposition was proved in \cite{Ku6} Propositions 
6.4.3, 6.4.5, 6.4.6.

\begin{Proposition} \label{KKP2}  
Suppose that $m \ell$ is in 
${\cal N}_{1}^{(N)}(K_{[\epsilon(m\ell)]})$. 
Then \\
{\rm (0)} $\kappa_{m, \ell} \in 
\Sel(O_{K}[1/m \ell], E[p^{N}])$.\\
{\rm (1)} $\partial_{r}(\kappa_{m,\ell})=
\phi_{r}(\kappa_{\frac{m}{r},\ell})$ for any prime divisor 
$r$ of $m$. \\
{\rm (2)} $\partial_{\ell}(\kappa_{m,\ell})=\delta_{m}t_{\ell,K}$. \\
{\rm (3)} Assume further that $m$ is admissible. 
Then $\phi_{r}(\kappa_{m, \ell})=0$ for any prime divisor 
$r$ of $m$. \\
{\rm (4)} Assume further that $m\ell$ is admissible, and 
$m \ell$ is in ${\cal N}_{1}^{(N)}(K_{[\epsilon(m\ell)+1]})$. 
Then we have
$$\phi_{\ell}(\kappa_{m, \ell})=-\delta_{m\ell} t_{\ell,K}.$$ 
\end{Proposition}

\section{Relations of Selmer groups} \label{section4}

In this section, we prove a generalized version of Theorem \ref{MMT1}.

\subsection{Injectivity theorem} \label{subsec41}

Suppose that $K$ is in ${\cal K}_{(p)}$ and 
that $m$ is in ${\cal N}_{1}^{(N)}(K)$.
For a prime divisor $r$ of $m$, we denote by 
$$w_{r}: {\cal H}_{r}^{2}(K) 
\longrightarrow \Sel(O_{K}, E[p^{N}])^{\vee}$$
the homomorphism which is the dual of 
$\Sel(O_{K}, E[p^{N}]) \longrightarrow \bigoplus_{v \mid r} E(K_{v}) 
\otimes {\bf Z}/p^{N}$.
Recall that ${\cal H}_{r}^{2}(K)$ is a free $R_{K}/p^{N}$-module of 
rank $1$, generated by $t_{r,K}$. 

\begin{Proposition} \label{GenerationTheorem}
We assume that 
$\delta_{m}$ is a unit of $R_{K}/p^{N}$ 
for some $m \in {\cal N}_{1}^{(N)}(K)$.  
Then the natural homomorphism 
$\oplus_{r \mid m} w_{r}: 
\bigoplus_{r \mid m}  {\cal H}_{r}^{2}(K)   
\longrightarrow \Sel(O_{K}, E[p^{N}])^{\vee}$
is surjective.
\end{Proposition}

\begin{Remark} \label{GenerationTheoremRem}
\begin{rm}
We note that $\delta_{m}$ is numerically computable, in principle. 
\end{rm}
\end{Remark}

\noindent Proof of Proposition \ref{GenerationTheorem}. 
Let $x$ be an arbitrary element in  
$\Sel(O_{K}, E[p^{N}])^{\vee}$. 
Let $w_{r}: {\cal H}_{r}^{2}(K)   
\longrightarrow \Sel(O_{K}, E[p^{N}])^{\vee}$ be 
the natural homomorphism 
for each $r \mid m$. 
We will prove that $x$ is in the submodule 
generated by all $w_{r}(t_{r,K})$ for $r \mid m$.
Using Lemma \ref{Lemma621}, we can take 
a prime $\ell \in ({\cal P}_{0}')^{(N)}(K(m)_{[1]}K_{n_{0}+N})$ 
such that 
$w_{\ell}(t_{\ell,K})=x$ and $\ell$ is prime to $m$. 
We consider the Kolyvagin derivative $\kappa_{m,\ell}$ 
which was defined in (\ref{KDer}).
Consider the exact sequence 
$$\Sel(O_{K}[1/m \ell], E[p^{N}])  
\stackrel{\partial}{\longrightarrow}  
\bigoplus_{\ell' \mid m \ell}  {\cal H}_{\ell'}^{2}(K)   
\stackrel{w_{K}}{\longrightarrow} \Sel(O_{K}, E[p^{N}])^{\vee}  
$$
(see Lemma \ref{FES0}) 
where $\partial=(\oplus \partial_{\ell'})_{\ell' \mid m \ell}$ and 
$w_{K}((z_{\ell'})_{\ell' \mid m \ell})=\sum_{\ell' \mid m \ell} 
w_{\ell'}(z_{\ell'})$.  
For each $r \mid m$ we define $\lambda_{r} \in R_{K}/p^{N}$ by 
$\partial_{r}(\kappa_{m,\ell})=
\lambda_{r}t_{r,K} \in {\cal H}_{r}^{2}(K)$.  
The above exact sequence and Proposition \ref{KKP1} (2) imply that   
$$\delta_{m} x +\sum_{r \mid m} \lambda_{r}
w_{r}(t_{r,K})=0$$
in $\Sel(O_{K}, E[p^{N}])^{\vee}$. 
Since we assumed that $\delta_{m}$ is a unit, $x$ is in the submodule 
generated by all $w_{r}(t_{r,K})$'s.
This completes the proof of Proposition \ref{GenerationTheorem}.

\vspace{5mm}

For a prime $\ell \in {\cal P}_{1}^{(N)}(K)$, 
we define 
$${\cal H}_{\ell,f}^{1}(K)=\bigoplus_{v \mid \ell} 
E(\kappa(v)) \otimes {\bf Z}/p^{N}.$$
Since $\kappa(v)={\bf F}_{\ell}$, 
$E(\kappa(v)) \otimes {\bf Z}/p^{N}$ is isomorphic to 
${\bf Z}/p^{N}$ and 
${\cal H}_{\ell,f}^{1}(K)$ is a free $R_{K}/p^{N}$-module of 
rank $1$. 

\begin{Corollary} \label{injectivitytheorem}
Suppose that $m=\ell_{1} \cdot ... \cdot \ell_{a}$ is 
in ${\cal N}_{1}^{(N)}(K)$. 
We assume that  
$\delta_{m}$ is a unit of $R_{K}/p^{N}$.  
Then the natural homomorphism 
$$s_{m}: \Sel(O_{K}, E[p^{N}]) 
\longrightarrow 
\bigoplus_{i=1}^{a} {\cal H}_{\ell_{i},f}^{1}(K)$$
is injective.
\end{Corollary}

\noindent Proof. This is obtained by taking the dual of 
the statement in Proposition \ref{GenerationTheorem}. 

\subsection{Relation matrices} \label{subsec42}

\begin{Theorem} \label{RelationTheorem}
Suppose that $m=\ell_{1} \cdot ... \cdot \ell_{a}$ is 
in ${\cal N}_{1}^{(N)}(K_{[a+1]})$. 
We assume that $m$ is admissible
and that 
$\delta_{m}$ is a unit of $R_{K}/p^{N}$.  
Then \\
{\rm (1)} $\Sel(O_{K}[1/m], E[p^{N}])$ is a free 
$R_{K}/p^{N}$-module of rank $a$. \\
{\rm (2)} $\{\kappa_{\frac{m}{\ell_{i}},\ell_{i}}\}_{1 \leq i \leq a}$ is 
a basis of $\Sel(O_{K}[1/m], E[p^{N}])$.  \\
{\rm (3)} The matrix 
\begin{equation} \label{relationmatrix}
{\cal A}=\left( \begin{array}{ccccc}
\delta_{\frac{m}{\ell_{1}}} & 
\phi_{\ell_{1}}(\kappa_{\frac{m}{\ell_{1}\ell_{2}},\ell_{2}}) &  & ... & 
\phi_{\ell_{1}}(\kappa_{\frac{m}{\ell_{1}\ell_{a}},\ell_{a}}) \\
\phi_{\ell_{2}}(\kappa_{\frac{m}{\ell_{1}\ell_{2}},\ell_{1}}) & 
\delta_{\frac{m}{\ell_{2}}} &  & ... & 
\phi_{\ell_{2}}(\kappa_{\frac{m}{\ell_{2}\ell_{a}},\ell_{a}}) \\
. &   &   & .   &   \\
. &   &   &  .  &   \\
\phi_{\ell_{a}}(\kappa_{\frac{m}{\ell_{1}\ell_{a}},\ell_{1}}) &
\phi_{\ell_{a}}(\kappa_{\frac{m}{\ell_{2}\ell_{a}},\ell_{2}}) &
   & ... & 
\delta_{\frac{m}{\ell_{a}}} 
\end{array} \right)
\end{equation}
is a relation matrix of $\Sel(E/K, E[p^{N}])^{\vee}$.
\end{Theorem}

In particular, if $a=2$, the above matrix is 
${\cal A}=\left( \begin{array}{cc}
\delta_{\ell_{2}} & 
\phi_{\ell_{1}}(g_{\ell_{2}}) \\
\phi_{\ell_{2}}(g_{\ell_{1}}) &
\delta_{\ell_{1}}
\end{array} \right)$. 
This is described in Remark 10.6 in \cite{Ku5} in the case of 
ideal class groups.

\vspace{5mm}

\noindent Proof of Theorem \ref{RelationTheorem} (1). 
By Proposition \ref{GenerationTheorem}, 
$\bigoplus_{i=1}^{a}  {\cal H}_{\ell_{i}}^{2}(K)   
\longrightarrow \Sel(O_{K}, E[p^{N}])^{\vee}$
is surjective.
Therefore, by Lemma \ref{FES0} we have an exact sequence
\begin{eqnarray} \label{ES3}
0 \longrightarrow \Sel(O_{K}, E[p^{N}])
&\longrightarrow &\Sel(O_{K}[1/m], E[p^{N}])
\stackrel{\partial}{\longrightarrow} 
\bigoplus_{i=1}^{a}  {\cal H}_{\ell_{i}}^{2}(K) \nonumber \\   
& \longrightarrow &\Sel(O_{K}, E[p^{N}])^{\vee}
\longrightarrow 0.
\end{eqnarray}
It follows that $\#\Sel(O_{K}[1/m], E[p^{N}])
=\#\bigoplus_{i=1}^{a}  {\cal H}_{\ell_{i}}^{2}(K)
=\#(R_{K}/p^{N})^{a}$. 

Let $m_{R_{K}}$ be the maximal ideal of $R_{K}$. 
By Lemma \ref{Lemma331}, 
$\Sel({\bf Z}[1/m], E[p^{N}]) 
\stackrel{\simeq}{\longrightarrow}
\Sel(O_{K}[1/m], E[p^{N}])^{\Gal(K/{\bf Q})}$ 
is bijective. 
Since $H^{0}({\bf Q}, E[p^{\infty}])=0$, the kernel of 
the multiplication by $p$ on $\Sel({\bf Z}[1/m], E[p^{N}])$ 
is $\Sel({\bf Z}[1/m], E[p])$. 
Therefore, we have an isomorphism
$\Sel(O_{K}[1/m], E[p^{N}])^{\vee} \otimes_{R_{K}}R_{K}/m_{R_{K}}
\simeq \Sel({\bf Z}[1/m], E[p])^{\vee}$. 
From the exact sequence 
$$0 \longrightarrow \Sel({\bf Z}, E[p])
\longrightarrow \Sel({\bf Z}[1/m], E[p])
\longrightarrow \bigoplus_{i=1}^{a} {\cal H}_{\ell_{i}}^{2}({\bf Q})
\longrightarrow \Sel({\bf Z}, E[p])^{\vee}
\longrightarrow 0,$$
and  ${\cal H}_{\ell_{i}}^{2}({\bf Q})=
H^{0}({\bf F}_{\ell_{i}}, E[p]) \simeq {\bf F}_{p}$,
we know that $\Sel({\bf Z}[1/m], E[p])$ is generated by $a$ elements. 
Therefore, by Nakayama's lemma, 
$\Sel(O_{K}[1/m], E[p^{N}])^{\vee}$  is generated by $a$ elements.
Since $\#\Sel(O_{K}[1/m], E[p^{N}])^{\vee}=\#(R_{K}/p^{N})^{a}$, 
$\Sel(O_{K}[1/m], E[p^{N}])^{\vee}$ is 
a free $R_{K}/p^{N}$-module of rank $a$. 
This shows that $\Sel(O_{K}[1/m], E[p^{N}])$ is also  
a free $R_{K}/p^{N}$-module of rank $a$ 
because $R_{K}/p^{N}$ is a Gorenstein ring. 

\vspace{5mm}

\noindent (2) 
We identify $\bigoplus_{i=1}^{a}  {\cal H}_{\ell_{i}}^{2}(K)$ 
with $(R_{K}/p^{N})^{a}$, using a basis 
$\{t_{\ell_{i},K}\}_{1 \leq i \leq a}$. 
Consider 
$\phi_{\ell_{i}}: \Sel(O_{K}[1/m], E[p^{N}])
\longrightarrow {\cal H}_{\ell_{i}}^{2}(K)$ 
and the direct sum of $\phi_{\ell_{i}}$, which we denote by $\Phi$;
$$\Phi=\oplus_{i=1}^{a} \phi_{\ell_{i}}: \Sel(O_{K}[1/m], E[p^{N}])
\longrightarrow \bigoplus_{i=1}^{a}{\cal H}_{\ell_{i}}^{2}(K)
\simeq (R_{K}/p^{N})^{a}.$$
Recall that 
$\kappa_{\frac{m}{\ell_{i}},\ell_{i}}$ is an element of  
$\Sel(O_{K}[1/m], E[p^{N}])$ (Proposition \ref{KKP2} (0)). 
By Proposition \ref{KKP2} (3), (4), we have 
$$\Phi(\kappa_{\frac{m}{\ell_{i}},\ell_{i}})=-\delta_{m} e_{i}$$
for each $i$ where $\{e_{i}\}_{1 \leq i \leq a}$ is the standard basis of 
the free module $(R_{K}/p^{N})^{a}$. 
Since we are assuming that $\delta_{m}$ is a unit, 
$\Phi$ is surjective. 
Since both the target and the source are free modules of the same rank, 
$\Phi$ is bijective. 
This implies Theorem \ref{RelationTheorem} (2). 

\vspace{5mm}

\noindent (3) 
Using the exact sequence (\ref{ES3}) and the isomorphism $\Phi$, 
we have an exact sequence 
$$(R_{n}/p^{N})^{a} \stackrel{\ \partial \circ \Phi^{-1}}
{\longrightarrow}
\bigoplus_{1 \leq i \leq a}  {\cal H}_{\ell_{i}}^{2}(K_{n})   
\stackrel{r}{\longrightarrow} \Sel(O_{K}, E[p^{N}])^{\vee}  
\longrightarrow  0.$$
We take a basis $\{-\delta_{m}e_{i}\}_{1 \leq i \leq a}$ of 
$(R_{n}/p^{N})^{a}$ and 
a basis $\{t_{\ell_{i},K}\}_{1 \leq i \leq a}$ of 
$\bigoplus_{1 \leq i \leq a} {\cal H}_{\ell_{i}}^{2}(K_{n})$. 
Then the $(i,j)$-component of the matrix corresponding to 
$\partial \circ \Phi^{-1}$ is 
$\partial_{\ell_{i}}(\kappa_{\frac{m}{\ell_{j}},\ell_{j}})$. 
If $i=j$, this is 
$\delta_{\frac{m}{\ell_{i}}}$ by Proposition \ref{KKP2} (2).   
If $i \neq j$, we have 
$\partial_{\ell_{i}}(\kappa_{\frac{m}{\ell_{j}},\ell_{j}})=
\phi_{\ell_{i}}(\kappa_{\frac{m}{\ell_{i}\ell_{j}},\ell_{j}})$ 
by Proposition \ref{KKP2} (1). 
This completes the proof of Theorem \ref{RelationTheorem}.

\begin{Remark}
\begin{rm}
Suppose that $\ell$ is in ${\cal P}_{1}^{(N)}(K)$. 
We define 
$$\Phi_{\ell}': H^{1}(K, E[p^{N}]) \longrightarrow 
{\cal H}_{\ell,f}^{1}(K)$$
as the composition of the natural map 
$H^{1}(K, E[p^{N}]) \longrightarrow 
\bigoplus_{v \mid \ell} H^{1}(K_{v}, E[p^{N}])$ 
and $\phi': 
H^{1}(K_{v}, E[p^{N}]) \longrightarrow 
H^{1}(\kappa(v), E[p^{N}])=E(\kappa(v)) \otimes {\bf Z}/p^{N}$ 
in (\ref{phiprime}).
For $m \in {\cal N}_{1}^{(N)}(K)$, we define 
$$\Phi_{m}': H^{1}(K, E[p^{N}]) \longrightarrow 
\bigoplus_{\ell \mid m}{\cal H}_{\ell,f}^{1}(K)$$
as the direct sum of $\Phi_{\ell}'$ for $\ell \mid m$. 
By definition, 
the restriction of $\Phi_{m}'$ to 
${\cal S}=\Sel(E/K, E[p^{N}])$ coincides with the canonical map 
$s_{m}$; 
\begin{equation}
(\Phi_{m}')_{\mid_{\cal S}}=s_{m}: 
\Sel(E/K, E[p^{N}]) \longrightarrow
\bigoplus_{\ell \mid m}{\cal H}_{\ell,f}^{1}(K) \ .
\end{equation}

Since ${\cal H}_{\ell,f}^{1}(K)$ and ${\cal H}_{\ell}^{2}(K)$ 
are Pontrjagin dual each other, we can take the dual basis 
$t^{{*}}_{\ell,K}$ of ${\cal H}_{\ell,f}^{1}(K)$ as an 
$R_{K}/p^{N}$-module from the basis $t_{\ell,K}$ of ${\cal H}_{\ell}^{2}(K)$. 
Under the assumptions of Theorem \ref{RelationTheorem},
using the basis $\{t^{{*}}_{\ell_{i},K}\}_{1 \leq i \leq a}$ 
of  $\bigoplus_{i=1}^{a}
{\cal H}_{\ell,f}^{1}(K)$, $\{t_{\ell_{i},K}\}_{1 \leq i \leq a}$ 
of  $\bigoplus_{i=1}^{a}
{\cal H}_{\ell_{i}}^{2}(K)$ and 
the isomorphism $\Phi_{m}'$, we have an exact sequence 
$\bigoplus_{\ell \mid m}{\cal H}_{\ell,f}^{1}(K) 
\stackrel{f}{\longrightarrow} 
\bigoplus_{i=1}^{a}
{\cal H}_{\ell_{i}}^{2}(K)
\longrightarrow \Sel(E/K, E[p^{N}])^{\vee} \longrightarrow 0$. 
Then the matrix corresponding to $f$ is an organizing matrix 
in the sense of Mazur and Rubin \cite{MR2} (cf. \cite{Ku6} \S 9). 
\end{rm} 
\end{Remark}

\section{Modified Kolyvagin systems and numerical examples} \label{section5}

\subsection{Modified Kolyvagin systems of Gauss sum type} \label{subsection51}

In \S \ref{subsec34} we constructed Kolyvagin systems 
$\kappa_{m,\ell}$ for $(m,\ell)$ such that $m\ell \in 
{\cal N}_{1}^{(N)}(K_{[\epsilon(m\ell)+1]})$. 
But the condition $\ell \in {\cal P}_{1}^{(N)}(K_{[\epsilon(m\ell)+1]})$ 
is too strict, and it is not suitable for numerical computation. 
In this subsection, we define a modified version of Kolyvagin systems of 
Gauss sum type for $(m, \ell)$ such that $m\ell \in 
{\cal N}_{1}^{(N)}(K)$. 

Suppose that $K$ is in ${\cal K}_{(p)}$. 
For each $\ell \in {\cal P}_{1}^{(N)}(K)$, we fix $t_{\ell} 
\in H^{0}({\bf F}_{\ell}, E[p^{N}])$ of order $p^{N}$, 
and consider $t_{\ell,K} \in {\cal H}^{2}_{\ell}(K)$, 
whose $\ell_{K}$-component is 
$t_{\ell} \otimes \zeta_{p^{N}}^{\otimes (-1)}$ and 
other components are zero. 
Using $t_{\ell,K}$, we regard $\partial_{\ell}$ and 
$\phi_{\ell}$ as homomorphisms 
$\partial_{\ell}: H^{1}(K, E[p^{N}]) \longrightarrow R_{K}/p^{N}$ and  
$\phi_{\ell}: H^{1}(K, E[p^{N}]) \longrightarrow R_{K}/p^{N}$.

We will define an element $\kappa_{m, \ell}^{q,q',z}$ in 
$\Sel(O_{K}[1/m \ell], E[p^{N}])$ for $(m, \ell)$ such that 
$m \ell \in {\cal N}_{1}(K)$ 
(and for some primes $q$, $q'$ and some $z$ in 
$\Sel(O_{K}[1/qq'], E[p^{N}])$). 
Consider $(m,\ell)$ such that $\ell$ is a prime and 
$m\ell \in {\cal N}_{1}(K)$. 
We take $n_{0}$ sufficiently large such that 
every prime of $K_{n_{0}}$ dividing $m\ell$ is inert in 
$K_{\infty}/K_{n_{0}}$. 
Then by Proposition \ref{KKP1} (1), for any 
$q \in ({\cal P}_{0}')^{(N)}(K(m\ell)_{[1]}K_{n_{0}+N})$, 
$\kappa_{m\ell, q} \in \Sel(O_{K}[1/m \ell q], E[p^{N}])$ satisfies 
$$\partial_{r}(\kappa_{m\ell, q})=
\phi_{r}(\kappa_{\frac{m \ell}{r},q})$$
for all $r$ dividing $m\ell$.
By Lemma \ref{Lemma621}, we can take $q$, 
$q' \in ({\cal P}_{0}')^{(N)}(K(m\ell)_{[1]}K_{n_{0}+N})$
satisfying \\ 
$\bullet$ $w_{K}(t_{q, K})=w_{K}(t_{q', K})$, and \\
$\bullet$ there is $z \in H^{1}_{f}(O_{K}[1/qq'], E[p^{N}])$ such that 
$\partial_{K}(z)=t_{q, K}-t_{q', K}$, 
$\phi_{\ell}(z)=1$ and $\phi_{r}(z)=0$ for any $r$ dividing $m$.

For any $m \in {\cal N}_{1}(K)$, 
let $\delta_{m}$ 
be the element defined in (\ref{Dcoeff}).  
We define 
\begin{equation} \label{1011}
\kappa_{m, \ell}^{q,q',z} = \kappa_{m\ell, q}-\kappa_{m\ell, q'} 
-\delta_{m\ell}z \ .
\end{equation}
By Proposition \ref{KKP1} (2), we have 
$\kappa_{m, \ell}^{q,q',z} \in \Sel(O_{K}[1/m \ell], E[p^{N}])$.

\begin{Proposition} \label{KKP3} 
{\rm (0)} $\kappa_{m, \ell}^{q,q',z}$ is in $\Sel(O_{K}[1/m \ell], E[p^{N}])$. \\
{\rm (1)} The element $\kappa_{m,\ell}^{q,q',z}$ satisfies 
$\partial_{r}(\kappa_{m,\ell}^{q,q',z})=
\phi_{r}(\kappa_{\frac{m}{r},\ell}^{q,q',z})$ 
for any prime divisor $r$ of $m$. \\
{\rm (2)} We further assume that $m \ell$ is admissible in the sense of 
the paragraph before Proposition \ref{KKP1}.  
Then we have $\phi_{r}(\kappa_{m, \ell}^{q,q',z})=0$ 
for any prime divisor $r$ of $m$. \\
{\rm (3)} Under the same assumptions as {\rm (2)}, 
$\phi_{\ell}(\kappa_{m, \ell}^{q,q',z})=-\delta_{m \ell}$ holds.
 \end{Proposition}

\noindent Proof. (1) Using the definition of 
$\kappa_{m,\ell}^{q,q',z}$ and Proposition \ref{KKP1} (1),
we have 
$\partial_{r}(\kappa_{m,\ell}^{q,q',z}) = 
\partial_{r}(\kappa_{m\ell, q}-\kappa_{m\ell, q'})=
\phi_{r}(\kappa_{\frac{m \ell}{r},q} -  \kappa_{\frac{m \ell}{r},q'})$.
Next, we use the definition of 
$\kappa_{\frac{m}{r},\ell}^{q,q',z}$ 
and $\phi_{r}(z)=0$ to get 
$\phi_{r}(\kappa_{\frac{m \ell}{r},q} -  \kappa_{\frac{m \ell}{r},q'})
=\phi_{r}(\kappa_{\frac{m}{r},\ell}^{q,q',z}+\delta_{\frac{m \ell}{r}}z)=
\phi_{r}(\kappa_{\frac{m}{r},\ell}^{q,q',z})$.
These computations imply (1). \\
(2) We have 
$\phi_{r}(\kappa_{m\ell, q})=\phi_{r}(\kappa_{m\ell, q'})=0$
by Proposition \ref{KKP1} (3). 
This together with $\phi_{r}(z)=0$ implies 
$\phi_{r}(\kappa_{m, \ell}^{q,q',z})=
\phi_{r}(\kappa_{m\ell, q}-\kappa_{m\ell, q'}-\delta_{m\ell}z)=0$. \\
(3) We again use Proposition \ref{KKP1} (3) to get 
$\phi_{\ell}(\kappa_{m\ell, q})=\phi_{\ell}(\kappa_{m\ell, q'})=0$.
Since $\phi_{\ell}(z)=1$, we have 
$\phi_{\ell}(\kappa_{m, \ell}^{q,q',z})=
\phi_{\ell}(\kappa_{m\ell, q}-\kappa_{m\ell, q'}-\delta_{m\ell}z)=
-\delta_{m\ell}$.
This completes the proof of Proposition \ref{KKP3}.

\subsection{Proof of Theorem \ref{MMT2}} \label{subsection52}

In this subsection we take $K={\bf Q}$. 
For $m \in {\cal N}^{(N)}={\cal N}^{(N)}({\bf Q})$, we consider 
$\delta_{m} \in {\bf Z}/p^{N}$, which is 
defined from $\vartheta_{{\bf Q}(m)}$ by 
(\ref{Dcoeff}). 
We define $\tilde{\delta}_{m} \in {\bf Z}/p^{N}$ by 
\begin{equation} \label{Dcoeff2}
\tilde{\theta}_{{\bf Q}(m)} \equiv \tilde{\delta}_{m}
\prod_{i=1}^{r}(\tau_{\ell_{i}}-1) \pmod{p^{N}, 
(\tau_{\ell_{1}}-1)^{2},...,(\tau_{\ell_{r}}-1)^{2}}
\end{equation}
where $m=\ell_{1}\cdot...\cdot \ell_{r}$. 
By (\ref{Res3}), $\tilde{\theta}_{{\bf Q}(m)} = u
\vartheta_{{\bf Q}(m)}$ for some unit $u \in R_{{\bf Q}(m)}^{\times}$. 
This together with (\ref{Dcoeff}) and (\ref{Dcoeff2}) implies that 
\begin{equation} \label{Dcoeff3}
\ord_{p}(\tilde{\delta}_{m})=\ord_{p}(\delta_{m}).
\end{equation}
We take a generator $\eta_{\ell} \in ({\bf Z}/\ell {\bf Z})^{\times}$
such that the image of $\sigma_{\eta_{\ell}}
\in \Gal({\bf Q}(\mu_{\ell})/{\bf Q}) \simeq ({\bf Z}/\ell)^{\times}$ in 
$\Gal({\bf Q}(\ell)/{\bf Q}) \simeq 
({\bf Z}/\ell)^{\times} \otimes {\bf Z}_{p}$ is 
$\tau_{\ell}$ which is the generator 
we took. 
Then, using (\ref{Dcoeff2}) and (\ref{ModularElement}), we can easily 
check that the equation (\ref{ModularSymbolDelta}) in \S 
\ref{subsec11} holds. 

\vspace{5mm}

In the rest of this subsection, we take $N=1$. 
We simply write ${\cal P}_{1}$ for ${\cal P}_{1}^{(1)}$, 
so 
$${\cal P}_{1}=\{\ell \in {\cal P}_{good} \mid 
\ell \equiv 1 \ \mbox{(mod $p$) and} \ 
E({\bf F}_{\ell}) \simeq {\bf Z}/p\}.$$
The set of squarefree products of 
primes in ${\cal P}_{1}$ is denoted by ${\cal N}_{1}$.

We first prove the following lemma which is related to 
the functional equation of an elliptic curve.  

\begin{Lemma} \label{10FunctionalEquation}
Let $\epsilon$ be the root number of $E$. 
Suppose that $m \in {\cal N}_{1}$ is $\delta$-minimal 
(for the definition of $\delta$-minimalness, see the 
paragraph before Conjecture \ref{Conj2}).  
Then we have $\epsilon=(-1)^{\epsilon(m)}$. 
\end{Lemma}

\noindent Proof. 
By the functional equation (1.6.2) in Mazur and Tate \cite{MT} 
and the above 
definition of $\tilde{\delta}_{m}$, we have 
$\epsilon(-1)^{\epsilon(m)} \tilde{\delta}_{m} \equiv \tilde{\delta}_{m}$ 
(mod $p$). 
Since $\tilde{\delta}_{m} \not \equiv 0$ (mod $p$) is equivalent to 
$\delta_{m} \not \equiv 0$ (mod $p$) by (\ref{Dcoeff3}), 
we get the conclusion. 

\vspace{5mm}

For each $\ell \in {\cal P}_{1}$, we fix a generator 
$t_{\ell} \in {\cal H}^{2}_{\ell}({\bf Q})=
H^{0}({\bf F}_{\ell}, E[p](-1)) \simeq {\bf Z}/p={\bf F}_{p}$, 
and regard $\phi_{\ell}$ as a map 
$\phi_{\ell}: H^{1}({\bf Q}, E[p]) \longrightarrow {\bf F}_{p}$. 
Note that the restriction of $\phi_{\ell}$ to $\Sel(E/{\bf Q}, E[p])$ 
is the zero map if and only if the natural map 
$s_{\ell}: \Sel(E/{\bf Q}, E[p]) \longrightarrow 
E({\bf F}_{\ell}) \otimes {\bf Z}/p \simeq {\bf F}_{p}$
is the zero map. 

\vspace{5mm}

\noindent I) Proof of Theorem \ref{MMT2} (1), (2).
 
\noindent Suppose that $\epsilon(m)=0$, namely $m=1$. 
Then $\delta_{1}=\theta_{\bf Q}$ mod $p$ 
$=L(E,1)/\Omega_{E}^{+}$ mod $p$. 
If $\delta_{1} \neq 0$, $\Sel(E/{\bf Q}, E[p])=0$ 
and $s_{1}$ is trivially bijective. 
Suppose next $\epsilon(m)=1$, so $m=\ell \in {\cal P}_{1}$.
It is sufficient to prove the next two propositions. 

\begin{Proposition} \label{10Prop1}
Assume that $\ell \in {\cal P}_{1}$ is $\delta$-minimal. 
Then $\Sel(E/{\bf Q}, E[p])$ is 
$1$-dimensional over ${\bf F}_{p}$, and 
$s_{\ell}: \Sel(E/{\bf Z}, E[p]) \longrightarrow {\bf F}_{p}$ 
is bijective. 
Moreover, the Selmer group 
$\Sel(E/{\bf Q}, E[p^{\infty}])^{\vee}$  
with respect to the $p$-power torsion points $E[p^{\infty}]$ 
is a free ${\bf Z}_{p}$-module of rank $1$, namely 
$\Sel(E/{\bf Q}, E[p^{\infty}])^{\vee} \simeq {\bf Z}_{p}$. 
\end{Proposition}

\noindent Proof. 
We first assume 
$\Sel(E/{\bf Q}, E[p])=0$ and will obtain the contradiction. 
We consider 
$\kappa_{1, \ell}^{q,q',z} = \kappa_{\ell, q}-\kappa_{\ell, q'} 
-\delta_{\ell}z$, which was defined in (\ref{1011}).  
By Proposition \ref{KKP1} (1), we know 
$\partial_{\ell}(\kappa_{1, \ell}^{q,q',z})=
\phi_{\ell}(g_{q}-g_{q'})$. 
Consider the exact sequence (see Lemma \ref{FES0})
$$0 \longrightarrow \Sel(E/{\bf Q}, E[p]) 
\longrightarrow \Sel({\bf Z}[1/r], E[p]) \longrightarrow 
{\cal H}^{2}_{r}({\bf Q})$$
for any $r \in {\cal P}_{1}$ where $\Sel({\bf Z}[1/r], E[p]) \longrightarrow 
{\cal H}^{2}_{r}({\bf Q}) \simeq {\bf F}_{p}$ is 
nothing but $\partial_{r}$. 
Since we assumed $\Sel(E/{\bf Q}, E[p])=0$, 
$\Sel({\bf Z}[1/r], E[p]) \longrightarrow 
{\cal H}^{2}_{r}({\bf Q}) \simeq {\bf F}_{p}$ 
is injective for any $r \in {\cal P}_{1}$. 
So 
$\partial_{q}(g_{q})=\delta_{1}=0$ implies 
that $g_{q}=0$. 
By the same method, we have $g_{q'}=0$. 
Therefore, $\partial_{\ell}(\kappa_{1, \ell}^{q,q',z})=
\phi_{\ell}(g_{q}-g_{q'})=0$, which implies 
that $\kappa_{1, \ell}^{q,q',z} \in \Sel(E/{\bf Q}, E[p])$. 

But Proposition \ref{KKP3} (3) tells us that 
$\phi_{\ell}(\kappa_{1, \ell}^{q,q',z})=-\delta_{\ell} \neq 0$. 
Therefore, $\kappa_{1, \ell}^{q,q',z} \neq 0$, which 
contradicts our assumption $\Sel(E/{\bf Q}, E[p]) = 0$.
Thus we get $\Sel(E/{\bf Q}, E[p]) \neq 0$. 

On the other hand, by Corollary \ref{injectivitytheorem}
we know that 
$s_{\ell}: \Sel(E/{\bf Q}, E[p]) \longrightarrow {\bf F}_{p}$ 
is injective, therefore bijective. 

By Lemma \ref{10FunctionalEquation}, the root 
number $\epsilon$ is $-1$. 
This shows that $\Sel(E/{\bf Q}, E[p^{\infty}])^{\vee}$ has 
positive ${\bf Z}_{p}$-rank by the parity conjecture proved by 
Nekov\'{a}\v{r} (\cite{Ne}). 
Therefore, we finally have 
$\Sel(E/{\bf Q}, E[p^{\infty}])^{\vee} \simeq {\bf Z}_{p}$, which 
completes the proof of Proposition \ref{10Prop1}. 

\vspace{5mm}

If we assume a slightly stronger condition on $\ell$, we 
also obtain the main conjecture. 
Let $\lambda'=\lambda^{an}$ be 
the analytic $\lambda$-invariant of the $p$-adic 
$L$-function $\vartheta_{{\bf Q}_{\infty}}$.
We put $n_{\lambda'} = \min\{n \in {\bf Z} \mid p^{n}-1 \geq \lambda'\}$.

\begin{Proposition} \label{10Prop2}
Suppose that there is $\ell \in {\cal P}_{1}$ such that 
$$\ell \equiv 1 \ \mbox{{\rm (mod $p^{n_{\lambda'}+2}$)} and} \ 
\tilde{\delta}_{\ell} \neq 0.$$  
Then the main conjecture for $(E, {\bf Q}_{\infty}/{\bf Q})$ is 
true and 
$\Sel(E/{\bf Q}_{\infty}, E[p^{\infty}])^{\vee}$
is generated by one element as a $\Lambda_{{\bf Q}_{\infty}}$-module. 
\end{Proposition}

\noindent Proof. 
We use our Euler system $g_{\ell}^{(K)}$ in \S \ref{subsec32}  
instead of $\kappa_{1, \ell}^{q,q',z}$ which was used 
in the proof of Proposition \ref{10Prop1}. 
Let $\lambda$ be the algebraic $\lambda$-invariant, namely the rank 
of $\Sel(E/{\bf Q}_{\infty}, E[p^{\infty}])^{\vee}$. 
Then $\lambda \leq \lambda'$ and 
$\vartheta_{{\bf Q}_{\infty}} \in 
\cc(\Sel(O_{{\bf Q}_{\infty}}, E[p^{\infty}])^{\vee})$
by Kato's theorem.

Put $K={\bf Q}_{n_{\lambda'}}$ and $f=p^{n_{\lambda'}}$. 
Consider the group ring $R_{K}/p={\bf F}_{p}[\Gal(K/{\bf Q})]$. 
We identify a generator $\gamma$ of $\Gal(K/{\bf Q})$ 
with $1+{\mathfrak t}$, and identify $R_{K}/p$ with 
${\bf F}_{p}[[{\mathfrak t}]]/({\mathfrak t}^{f})$.
The norm $N_{\Gal(K/{\bf Q})}=
\Sigma_{i=0}^{f-1} 
\gamma^{i}$ is ${\mathfrak t}^{f-1}$ 
by this identification, so 
our assumption $\lambda' \leq f-1$ implies that 
the corestriction map 
$\Sel(E/K, E[p]) \longrightarrow 
\Sel(E/{\bf Q}, E[p])$ is the zero map
because $\lambda \leq \lambda'$. 
Therefore, we have ${\bf Q}_{[1]} \subset K$.
Since $p^{n_{\lambda'}+1}-p^{n_{\lambda'}} > 
p^{n_{\lambda'}}-1 \geq \lambda' \geq \lambda$, 
the corestriction map $\Sel(E/{\bf Q}_{n_{\lambda'}+1}, E[p]) 
\longrightarrow 
\Sel(E/{\bf Q}_{n_{\lambda'}}, E[p])=\Sel(E/K, E[p])$ 
is also the zero map. 
This shows that ${\bf Q}_{[2]} \subset {\bf Q}_{n_{\lambda'}+1}$.

Our assumption 
$\ell \equiv 1$ (mod $p^{n_{\lambda'}+2}$) implies that 
$\ell$ splits completely in ${\bf Q}_{n_{\lambda'}+1}$, so 
we have $\ell \in {\cal P}_{1}({\bf Q}_{[2]})
={\cal P}_{1}(K_{[1]})$. 
Therefore, we can define 
$$g_{\ell}^{(K)} \in 
\Sel(O_{K}[1/\ell], E[p])$$ 
in \S \ref{subsec32}. 
Since $\ell \in {\cal P}_{1}({\bf Q}_{[2]})$, we also have 
$$\phi_{\ell}(g_{\ell}^{({\bf Q})})=-\delta_{\ell}^{({\bf Q})}=-\delta_{\ell}$$ 
by Proposition \ref{KKP2} (4).  
It follows from our assumption $\delta_{\ell} \neq 0$ that 
$g_{\ell}^{({\bf Q})} \neq 0$. 
Since $\Cor_{K/{\bf Q}}(g_{\ell}^{(K)})=g_{\ell}^{({\bf Q})}$ 
and the natural map 
$i: \Sel({\bf Z}[1/\ell], E[p]) \longrightarrow \Sel(O_{K}[1/\ell], E[p])$ 
is injective, we get
$$i(g_{\ell}^{({\bf Q})})=N_{\Gal(K/{\bf Q})} g_{\ell}^{(K)}
=
{\mathfrak t}^{f-1}g_{\ell}^{(K)}
\neq 0.$$ 

Consider $\partial_{\ell}: \Sel(O_{K}[1/\ell], E[p]) 
\longrightarrow R_{K}/p$. 
By definition, we have 
$\partial_{\ell}(g_{\ell}^{(K)}) = u{\mathfrak t}^{\lambda'}$ 
for some unit $u$ of $R_{K}/p$. 
This shows that $\partial_{\ell}({\mathfrak t}^{f-\lambda'}g_{\ell}^{(K)}) = 0$, 
which implies that ${\mathfrak t}^{f-\lambda'}g_{\ell}^{(K)} \in 
\Sel(E/K, E[p])$.  
The fact ${\mathfrak t}^{f-1}g_{\ell}^{(K)} \neq 0$ implies 
the submodule generated by ${\mathfrak t}^{f-\lambda'}g_{\ell}^{(K)}$ 
is isomorphic to $R_{K}/(p, {\mathfrak t}^{\lambda'})$ as an $R_{K}$-module.
Namely, we have 
$$\Sel(E/K, E[p]) \supset 
\langle {\mathfrak t}^{f-\lambda'}g_{\ell}^{(K)} \rangle 
\simeq R_{K}/(p, {\mathfrak t}^{\lambda'}).$$ 
This implies that $\lambda=\lambda'$, and 
$\Sel(E/K, E[p]) \simeq R_{K}/(p, {\mathfrak t}^{\lambda})$. 
Therefore, we have 
$\Sel(E/{\bf Q}_{\infty}, E[p])^{\vee} 
\simeq \Lambda_{{\bf Q}_{\infty}}/(p, \vartheta_{{\bf Q}_{\infty}})$. 
This together with Kato's theorem we mentioned implies that  
$\Sel(E/{\bf Q}_{\infty}, E[p^{\infty}])^{\vee}
\simeq \Lambda_{{\bf Q}_{\infty}}/(\vartheta_{{\bf Q}_{\infty}})$.

\vspace{5mm}

\noindent II) Proof of Theorem \ref{MMT2} (3). 

\noindent Suppose that $m =\ell_{1}\ell_{2} \in {\cal N}_{1}$ 
and $m$ is $\delta$-minimal. 
As in the proof of Proposition \ref{10Prop1}, we 
assume $\Sel(E/{\bf Q}, E[p])=0$ and will get the contradiction.  
We consider 
$\kappa_{\ell_{1}, \ell_{2}}^{q,q',z}$ defined in (\ref{1011}). 
Consider the exact sequence (see Lemma \ref{FES0})
$$0 \longrightarrow \Sel(E/{\bf Q}, E[p]) 
\longrightarrow \Sel({\bf Z}[1/\ell_{1}\ell_{2} q q'], E[p]) 
\stackrel{\partial}{\longrightarrow} 
\bigoplus_{v \in \{\ell_{1}, \ell_{2},q,q'\}} {\cal H}^{2}_{v}({\bf Q}).$$
By the same method as the proof of Proposition \ref{10Prop1}, 
$g_{q}=g_{q'}=0$. 
Therefore, $\partial_{\ell_{1}}(\kappa_{\ell_{1},q}-\kappa_{\ell_{1},q'})=
\phi_{\ell_{1}}(g_{q}-g_{q'})=0$ 
by Proposition \ref{KKP1} (1). 
We have  
$\partial_{q}(\kappa_{\ell_{1},q})=\delta_{\ell_{1}}=0$,  
$\partial_{q}(\kappa_{\ell_{1},q'})=0$, 
$\partial_{q'}(\kappa_{\ell_{1},q})=0$, 
$\partial_{q'}(\kappa_{\ell_{1},q'})=\delta_{\ell_{1}}=0$.
Therefore,   
$\partial(\kappa_{\ell_{1},q}-\kappa_{\ell_{1},q'})=0$. 
This together with $\Sel(E/{\bf Q}, E[p])=0$ shows that 
$\kappa_{\ell_{1},q}-\kappa_{\ell_{1},q'}=0$.
Therefore, using Proposition \ref{KKP1} (1), we have 
$$\partial_{\ell_{2}}(\kappa_{\ell_{1}, \ell_{2}}^{q,q',z})
=\partial_{\ell_{2}}(\kappa_{m, q}-\kappa_{m, q'})
=\phi_{\ell_{2}}(\kappa_{\ell_{1}, q}-\kappa_{\ell_{1}, q'})=0.$$
By the same method as the above proof of 
$\kappa_{\ell_{1},q}-\kappa_{\ell_{1},q'}=0$, 
we get $\kappa_{1, \ell_{2}}^{q,q',z}=0$.
This implies that $\partial_{\ell_{1}}(\kappa_{\ell_{1}, \ell_{2}}^{q,q',z})
=\phi_{\ell_{1}}(\kappa_{1, \ell_{2}}^{q,q',z})=0$
by Proposition \ref{KKP3} (1). 
It follows that $\partial(\kappa_{\ell_{1}, \ell_{2}}^{q,q',z})=0$, which 
implies $\kappa_{\ell_{1}, \ell_{2}}^{q,q',z} \in \Sel(E/{\bf Q}, E[p])$. 
But this is a contradiction because 
we assumed $\Sel(E/{\bf Q}, E[p])=0$ and 
$$\phi_{\ell_{2}}(\kappa_{\ell_{1}, \ell_{2}}^{q,q',z})=-\delta_{m} \neq 0$$
by Proposition \ref{KKP3} (3). 
Thus, we get $\Sel(E/{\bf Q}, E[p]) \neq 0$.

Now the root number is $1$ by 
Lemma \ref{10FunctionalEquation}, 
therefore, by the parity conjecture proved by 
Nekov\'{a}\v{r} (\cite{Ne}), we obtain 
$\dim_{{\bf F}_{p}}\Sel(E/{\bf Q}, E[p]) \geq 2$. 
On the other hand, 
by Corollary \ref{injectivitytheorem} we know that  
$s_{m}: \Sel(E/{\bf Q}, E[p]) \longrightarrow ({\bf F}_{p})^{\oplus 2}$
is injective. 
Therefore, the injectivity of $s_{m}$ 
implies the bijectivity of $s_{m}$. 
This completes the proof of Theorem \ref{MMT2} (3).

\vspace{5mm}

We give a simple corollary. 

\begin{Corollary} \label{10Prop3C}
Suppose that there is $m \in {\cal N}_{1}$ such that 
$m$ is $\delta$-minimal and $\epsilon(m)=2$. 
We further assume that 
the analytic $\lambda$-invariant $\lambda'$ is $2$. 
Then the main conjecture for $(E, {\bf Q}_{\infty}/{\bf Q})$ holds.  
\end{Corollary}

\noindent Proof. 
Put ${\mathfrak t}=\gamma-1$ and identify 
$\Lambda_{{\bf Q}_{\infty}}/p$ with 
${\bf F}_{p}[[{\mathfrak t}]]$. 
Let ${\cal A}$ be the relation matrix of 
$S=\Sel(E/{\bf Q}_{\infty}, E[p^{\infty}])^{\vee}$. 
Since $S/(p,{\mathfrak t})=\Sel(E/{\bf Q}, E[p])^{\vee} \simeq
{\bf F}_{p} \oplus {\bf F}_{p}$, 
${\mathfrak t}^{2}$ divides $\det {\cal A}$ mod $p$.
Therefore, the algebraic $\lambda$-invariant is also $2$. 
This implies the main conjecture because 
$\det {\cal A}$ divides $\vartheta_{{\bf Q}_{\infty}}$ in 
$\Lambda_{{\bf Q}_{\infty}}$ (Kato \cite{KK1}). 

\vspace{5mm}

\noindent III) Proof of Theorem \ref{MMT2} (4).
 
\begin{Lemma} \label{10L31}
Suppose that $\ell$, $\ell_{1}$, $\ell_{2}$ are 
distinct primes in ${\cal P}_{1}$ satisfying 
$\delta_{\ell}=\delta_{\ell \ell_{1}}=\delta_{\ell \ell_{2}}=0$. 
Assume also that  
$s_{\ell}: \Sel(E/{\bf Q}, E[p]) \longrightarrow {\bf F}_{p}$ 
is bijective, and that $\ell \ell_{1}$, $\ell \ell_{2}$ are both admissible.
We take $q$, $q'$ such that they satisfy the conditions when we defined 
$\kappa_{\ell_{1}\ell_{2},\ell}^{q,q',z}$.
Then we have\\
{\rm (1)} $\Sel(E/{\bf Q}, E[p]) = \Sel({\bf Z}[1/\ell], E[p])$, \\
{\rm (2)} $\kappa_{\ell_{1},\ell}^{q,q',z}=0$, 
$\kappa_{\ell_{2},\ell}^{q,q',z}=0$, and \\
{\rm (3)} $\kappa_{\ell_{1}\ell_{2},\ell}^{q,q',z} \in \Sel(E/{\bf Q}, E[p])$.
\end{Lemma}

\noindent Proof. (1) Since $s_{\ell}$ is bijective, taking the dual, we get 
the bijectivity of ${\cal H}^{2}_{\ell}({\bf Q})
\longrightarrow \Sel(E/{\bf Q}, E[p])^{\vee}=\Sel({\bf Z}, E[p])^{\vee}$. 
By the exact sequence
$$0 \longrightarrow \Sel({\bf Z}, E[p]) 
\longrightarrow \Sel({\bf Z}[1/\ell], E[p]) 
\stackrel{\partial_{\ell}}{\longrightarrow} 
{\cal H}^{2}_{\ell}({\bf Q})
\longrightarrow \Sel({\bf Z}, E[p])^{\vee}
\longrightarrow 0$$
 in Lemma \ref{FES0}, 
we get $\Sel(E/{\bf Q}, E[p])=
\Sel({\bf Z}, E[p]) = \Sel({\bf Z}[1/\ell], E[p])$. 

\noindent (2) We first note that the bijectivity of 
$s_{\ell}: \Sel(E/{\bf Q}, E[p]) \longrightarrow {\bf F}_{p}$ 
implies the bijectivity of 
$\phi_{\ell}: \Sel(E/{\bf Q}, E[p]) \longrightarrow {\bf F}_{p}$. 
Since $\partial_{q}(\kappa_{\ell, q})=\delta_{\ell}=0$, 
$\kappa_{\ell,q} \in  \Sel({\bf Z}[1/\ell], E[p])=
\Sel(E/{\bf Q}, E[p])$
where we used the property (1) which we have just proved. 
Proposition \ref{KKP1} (3) implies 
$\phi_{\ell}(\kappa_{\ell,q})=0$, which 
implies $\kappa_{\ell,q}=0$ by the bijectivity of $\phi_{\ell}$. 
By the same method, we have $\kappa_{\ell,q'}=0$. 
Therefore, we have 
$$\kappa_{1,\ell}^{q,q',z}=
\kappa_{\ell,q}-\kappa_{\ell,q'}-\delta_{\ell} z=0.$$
Therefore, Proposition \ref{KKP3} (1) implies 
$\partial_{\ell_{1}}(\kappa_{\ell_{1},\ell}^{q,q',z})=
\phi_{\ell_{1}}(\kappa_{1,\ell}^{q,q',z})=0$.
This implies 
$\kappa_{\ell_{1},\ell}^{q,q',z} \in  \Sel({\bf Z}[1/\ell], E[p])=
\Sel(E/{\bf Q}, E[p])$. 
Using Proposition \ref{KKP3} (3), we have 
$$\phi_{\ell}(\kappa_{\ell_{1},\ell}^{q,q',z})=-\delta_{\ell \ell_{1}}=0,$$
which implies $\kappa_{\ell_{1},\ell}^{q,q',z}=0$ by the 
bijectivity of $\phi_{\ell}$. 
The same proof works for $\kappa_{\ell_{2},\ell}^{q,q',z}$.

\noindent (3) It follows from Proposition \ref{KKP3} (1) and Lemma \ref{10L31} (2) 
that $\partial_{\ell_{i}}(\kappa_{\ell_{1}\ell_{2},\ell}^{q,q',z})=
\phi_{\ell_{i}}(\kappa_{\frac{\ell_{1}\ell_{2}}{\ell_{i}},\ell}^{q,q',z})=
0$ for each $i=1$, $2$. 
This implies $\kappa_{\ell_{1}\ell_{2},\ell}^{q,q',z} \in \Sel({\bf Z}[1/\ell], E[p])$. 
Using $\Sel({\bf Z}[1/\ell], E[p])=
\Sel(E/{\bf Q}, E[p])$ which we proved in (1), 
we get the conclusion. 
This completes the proof of Lemma \ref{10L31}.

\vspace{5mm}

We next prove Theorem \ref{MMT2} (4).
Assume that $m =\ell_{1}\ell_{2} \ell_{3} \in {\cal N}_{1}$, 
$m$ is $\delta$-minimal, $m$ is admissible, and 
$s_{\ell_{i}}: \Sel(E/{\bf Q}, E[p]) \longrightarrow {\bf F}_{p}$
is surjective for each $i=1$, $2$. 

We assume $\dim_{{\bf F}_{p}} \Sel(E/{\bf Q}, E[p]) =1$ and 
will get the contradiction.  
By this assumption, $s_{\ell_{i}}: \Sel(E/{\bf Q}, E[p]) 
\longrightarrow {\bf F}_{p}$ for each $i=1$, $2$ is bijective. 
This implies that 
$\phi_{\ell_{i}}: \Sel(E/{\bf Q}, E[p]) 
\longrightarrow {\bf F}_{p}$ for each $i=1$, $2$ is also bijective. 
By Lemma \ref{10L31} (3) we get 
$\kappa_{\ell_{2}\ell_{3},\ell_{1}}^{q,q',z} 
\in \Sel(E/{\bf Q}, E[p])$, 
taking $q$, $q'$ satisfying the conditions when we defined this element. 
By Proposition \ref{KKP3} (3), we have 
$\phi_{\ell_{1}}(\kappa_{\ell_{2}\ell_{3},\ell_{1}}^{q,q',z})=-\delta_{m} \neq 0$, 
which implies $\kappa_{\ell_{2}\ell_{3},\ell_{1}}^{q,q',z} \neq 0$. 
But by Proposition \ref{KKP3} (2), we have 
$\phi_{\ell_{2}}(\kappa_{\ell_{2}\ell_{3},\ell_{1}}^{q,q',z})=0$. 
This contradicts the bijectivity of $\phi_{\ell_{2}}$. 
Therefore, we obtain 
$\dim_{{\bf F}_{p}} \Sel(E/{\bf Q}, E[p]) > 1$. 

By Lemma \ref{10FunctionalEquation} and our assumption 
that $m$ is $\delta$-minimal, 
we know that the root number $\epsilon$ is $-1$. 
This shows that $\dim_{{\bf F}_{p}} \Sel(E/{\bf Q}, E[p]) \geq 3$
by the parity conjecture proved by Nekov\'{a}\v{r} (\cite{Ne}). 
On the other hand, Corollary \ref{injectivitytheorem} implies that 
$\dim_{{\bf F}_{p}} \Sel(E/{\bf Q}, E[p]) \leq 3$ 
and  
$s_{m}: 
\Sel(E/{\bf Q}, E[p]) \longrightarrow {\bf F}_{p}^{\oplus 3}$ 
is injective. 
Therefore,   
the above map $s_{m}$ is bijective. 
This completes the proof of Theorem \ref{MMT2} (4).

\subsection{Numerical examples} \label{subsection53}

In this section, we give several numerical examples.  

Let $E=X_{0}(11)^{(d)}$ be the quadratic twist of $X_{0}(11)$ by $d$, 
namely $dy^{2}=x^{3}-4x^{2}-160x-1264$. 
We take $p=3$. 
Then if $d \equiv 1$ (mod $p$), $p$ is a good ordinary prime which is 
not anomalous (namely $a_{p}(=a_{3})$ for $E$ satisfies 
$a_{p} \not \equiv 1$ (mod $p$)), 
and $p=3$ does not divide $\Tam(E)$, and the Galois representation on 
$T_{3}(E)$ is surjective.  
In the following examples, we checked $\mu'=0$ 
where $\mu'$ is the analytic $\mu$-invariant. 
Then this implies that the algebraic $\mu$-invariant 
is also zero (Kato \cite{KK1} Theorem 17.4 (3)) 
under our assumptions. 
In the computations of $\tilde{\delta}_{m}$ below, 
we have to fix a generator of 
$\Gal({\bf Q}(\ell)/{\bf Q}) \simeq ({\bf Z}/\ell{\bf Z})^{\times}$ 
for a prime $\ell$. 
We always take the least primitive root $\eta_{\ell}$ 
of $({\bf Z}/\ell{\bf Z})^{\times}$. 
We compute $\tilde{\delta}_{m}$ using the formula in 
(\ref{ModularSymbolDelta}). 

\vspace{5mm}

\noindent {\bf (1)} $d=13$.  
We take $N=1$. 
Since $\tilde{\delta}_{7}=20 \not \equiv 0$ (mod $3$), we know that 
$\Fitt_{1, {\bf F}_{3}}(\Sel(E/{\bf Q}, E[3])^{\vee})={\bf F}_{3}$ 
by Theorem \ref{HigherFittingIdeal}, so 
$\Sel(E/{\bf Q}, E[3])$
is generated by one element. 

The root number is $\epsilon=(\frac{13}{11})=-1$, so 
$L(E,1)=0$. 
We compute ${\cal P}_{1}=\{7, 31, 73,...\}$. 
Therefore, $\tilde{\delta}_{7} \not \equiv 0$ (mod $3$) implies 
$\Sel(E/{\bf Q}, E[3]) \simeq {\bf F}_{3}$ and 
$$\Sel(E/{\bf Q}, E[3^{\infty}])^{\vee} \simeq {\bf Z}_{3}$$
by Proposition \ref{10Prop1}.  
Also, it is easily computed that $\lambda'=1$ in this case. 
This implies that  
$\Sel(E/{\bf Q}_{\infty}, E[3^{\infty}])^{\vee} \simeq {\bf Z}_{3}$, so 
the main conjecture also holds. 

We can find a point $P=(7045/36 , -574201/216)$ of infinite order  on 
the minimal Weierstrass model 
$ y^2 + y = x^3 - x^2 - 1746x - 50295$ of $E=X_{0}(11)^{(13)}$. 
Therefore, we know $\sha(E/{\bf Q})[3^{\infty}]=0$.  
We can easily check that $E({\bf F}_{7})$ is cyclic of order $6$, and 
that the image of the point $P$  
in $E({\bf F}_{7})/3E({\bf F}_{7})$ is non-zero. 
So we also checked numerically that $s_{7}: \Sel(E/{\bf Q}, E[3]) 
\longrightarrow E({\bf F}_{7})/3E({\bf F}_{7})$ is  
bijective as Proposition \ref{10Prop1} claims.  

\vspace{5mm}

\noindent {\bf (2)} $d=40$. 
We know $\epsilon=(\frac{40}{11})=-1$. 
We take $N=1$.
We can compute ${\cal P}_{1}=\{7, 67, 73,...\}$, and 
$\tilde{\delta}_{7}=-40 \not \equiv 0$ (mod $3$). 
This implies that  
$\Sel(E/{\bf Q}, E[3]) \simeq {\bf F}_{3}$ and 
$\Sel(E/{\bf Q}, E[3^{\infty}])^{\vee} \simeq {\bf Z}_{3}$
by Proposition \ref{10Prop1}.  
 
In this case, we know $\lambda'=7$. 
Therefore, $n_{\lambda'}=2$. 
We can check $5347 \in {\cal P}_{1}$ 
(where $5347 \equiv 1$ (mod $3^{5}$)) 
and 
$\tilde{\delta}_{5347}=-412820 \not \equiv 0$ (mod $3$). 
Therefore, the main conjecture holds by 
Proposition \ref{10Prop2}.  
In this case, we can check that the $p$-adic $L$-function 
$\vartheta_{{\bf Q}_{\infty}}$ is 
divisible by $(1+{\mathfrak t})^{3}-1$, so we have 
$$\rank_{{\bf Z}_{3}} \Sel(E/{\bf Q}_{1}, E[3^{\infty}])^{\vee} =3$$ 
where ${\bf Q}_{1}$ is the first layer of ${\bf Q}_{\infty}/{\bf Q}$.

\vspace{5mm}

In the following, for a prime $\ell \in {\cal P}$, 
we take a generator $\tau_{\ell}$ of $\Gal({\bf Q}(\ell)/{\bf Q}) \simeq 
({\bf Z}/\ell {\bf Z})^{\times}$ and 
put $S=\tau_{\ell}-1$. 
We write 
$\vartheta_{{\bf Q}(\ell)}=\Sigma a_{i}^{(\ell)} S^{i}$ 
where $a_{i}^{(\ell)} \in {\bf Z}_{p}$.
Note that $\tilde{\delta}_{\ell}=a_{1}^{(\ell)}$. 

\vspace{5mm}

\noindent {\bf (3)} $d=157$. 
We know $\epsilon=(\frac{157}{11})=1$ and 
$L(E,1)/\Omega_{E}^{+}=45$. 
We take $N=1$.
We compute $a_{2}^{(37)} = -14065/2 \not \equiv 0$ (mod $3$). 
Since $37 \equiv 1$ (mod $3^{2}$), $c_{2}=2-1=1$ and 
$a_{2}^{(37)}$ is in 
$\Fitt_{2, {\bf F}_{3}}(\Sel(E/{\bf Q}, E[3])^{\vee})$ 
by Theorem \ref{HigherFittingIdeal}, which 
implies that 
$\Fitt_{2, {\bf F}_{3}}(\Sel(E/{\bf Q}, E[3]))={\bf F}_{3}$. 
Therefore, $\Sel(E/{\bf Q}, E[3])$
is generated by at most two elements.  

We compute 
${\cal P}_{1}=\{7, 67, 73, 127,...\}$. 
Since $127 \equiv 1$ (mod $7$), $7 \times 127$ is admissible. 
We compute $\tilde{\delta}_{7 \times 127}=83165 \not \equiv 0$ (mod $3$). 
Therefore, $7 \times 127$ is $\delta$-minimal. 
It follows from Theorem \ref{MMT2} (3) that 
$\Sel(E/{\bf Q}, E[3]) \simeq {\bf F}_{3} \oplus {\bf F}_{3}$.
In this example, we can check $\lambda'=2$, so Corollary \ref{10Prop3C} together with 
the above computation implies the 
main conjecture. 
Since $L(E,1)/\Omega_{E}^{+}=45 \neq 0$, $\rank E({\bf Q})=0$ by 
Kato, which implies $\Sel(E/{\bf Q}, E[3^{\infty}])=\sha(E/{\bf Q})[3^{\infty}]$. 
Since $45 \in \Fitt_{0, {\bf Z}_{3}}(\Sel(E/{\bf Q}, E[3^{\infty}])^{\vee})$, 
we have $\#\sha(E/{\bf Q})[3^{\infty}] \leq 9$, and 
$$\sha(E/{\bf Q})[3^{\infty}] \simeq {\bf Z}/3{\bf Z} \oplus {\bf Z}/3{\bf Z}.$$

\vspace{5mm}

\noindent {\bf (4)} $d=265$. 
In this case, $\epsilon=(\frac{265}{11})=1$ and $L(E,1)=0$. 
We take $N=1$.
As in Example (3), we compute 
$a_{2}^{(37)} = 16985 \not \equiv 0$ (mod $3$), which  
implies that $\Sel(E/{\bf Q}, E[3])$
is generated by at most two elements as above.  
We compute ${\cal P}_{1}=\{7, 13, 31, 67, 103, 109, 127,...\}$. 
For an admissible pair $\{7, 127\}$, we have 
$\tilde{\delta}_{7 \times 127}=-138880 \not \equiv 0$ (mod $3$). 
Therefore, $7 \times 127$ is $\delta$-minimal and 
$\Sel(E/{\bf Q}, E[3]) \simeq {\bf F}_{3} \oplus {\bf F}_{3}$
by Theorem \ref{MMT2} (3).
Since $\lambda'=2$ in this case, by Corollary \ref{10Prop3C} we know that 
the main conjecture holds.
 
Since $L(E,1)=0$, we know 
$\rank \Sel(E/{\bf Q}, E[3^{\infty}])^{\vee} > 0$ by the main conjecture. 
This implies that 
$$\Sel(E/{\bf Q}, E[3^{\infty}])^{\vee} \simeq {\bf Z}_{3} \oplus {\bf Z}_{3}.$$ 

Now $E$ has a minimal Weierstrass model 
$y^{2}+y=x^{3}-x^{2}-725658x-430708782$. 
We can find rational points 
$P=(2403,108146)$ and $Q=(5901,-448036)$ on this curve. 
We can also easily check that $E({\bf F}_{7})$ is cyclic group of order $6$ and 
$E({\bf F}_{31})$ is cyclic of order $39$.
The image of $P$ in $E({\bf F}_{7})/3E({\bf F}_{7}) \simeq {\bf Z}/3{\bf Z}$ 
is $0$ (the identity element), and the image of $Q$ in 
$E({\bf F}_{7})/3E({\bf F}_{7})\simeq {\bf Z}/3{\bf Z}$ is of order $3$. 
On the other hand,   
the images of $P$ and $Q$ in 
$E({\bf F}_{31})/3E({\bf F}_{31})\simeq {\bf Z}/3{\bf Z}$ 
do not vanish and coincide. 
This shows that $P$ and $Q$ are linearly independent over ${\bf Z}_{3}$. 
Therefore, 
$$\rank E({\bf Q})=2 \ \mbox{and} \ 
\sha(E/{\bf Q})[3^{\infty}]=0.$$  

In the above argument we considered the images 
of $E({\bf Q})$ in $E({\bf F}_{7})/3E({\bf F}_{7})$ 
and $E({\bf F}_{31})/3E({\bf F}_{31})$. 
What we explained above implies that the natural map 
$s_{7 \times 31}: E({\bf Q})/3E({\bf Q}) 
\longrightarrow E({\bf F}_{7})/3E({\bf F}_{7}) \oplus 
E({\bf F}_{31})/3E({\bf F}_{31})$ is bijective.  
In this example, $\tilde{\delta}_{7 \times 31}=-15290 \not \equiv 0$ (mod $3$), so Conjecture \ref{Conj2} holds for $m=7 \times 31$.

 
\vspace{5mm}

\noindent {\bf (5)} $d=853$. 
We know $\epsilon=(\frac{853}{11})=-1$. 
Take $N=1$ at first. 
For $\ell=271$, we have 
$a_{3}^{(271)} = 900852395/2 \not \equiv 0$ (mod $3$), which 
implies that $\dim_{{\bf F}_{3}} \Sel(E/{\bf Q}, E[3]) \leq 3$. 
We compute 
${\cal P}_{1}=\{7, 13, 67, 103, 109,...,463,...\}$.
We can find a rational point 
$P=(1194979057/51984, 40988136480065/11852352)$ 
on the minimal Weierstrass equation 
$y^{2} + y = x^{3} - x^{2} - 7518626x - 14370149745$ of 
$E=X_{0}(11)^{(853)}$. 
We know that $E({\bf F}_{7})$ is cyclic of order $6$, and 
$E({\bf F}_{13})$ is cyclic of order $18$.
Both of the images of $P$ in $E({\bf F}_{7})/3E({\bf F}_{7})$ and 
$E({\bf F}_{13})/3E({\bf F}_{13})$ are of order $3$.
Therefore, $s_{\ell}: \Sel(E/{\bf Q}, E[3]) \longrightarrow 
E({\bf F}_{\ell})/3E({\bf F}_{\ell})$ is surjective for each 
$\ell=7$, $13$. 
Since $13=-1 \in ({\bf F}_{7}^{\times})^{3}$, 
$463=1 \in ({\bf F}_{7}^{\times})^{3}$ and 
$463=8 \in ({\bf F}_{13}^{\times})^{3}$, 
$\{7,13,463\}$ is admissible. 
We can compute 
$\tilde{\delta}_{7 \times 13 \times 463}=-8676400 
\not \equiv 0$ (mod $3$), and can check that  
$m=7 \times 13 \times 463$ is $\delta$-minimal. 
By Theorem \ref{MMT2} (4), we have 
\begin{equation} \label{NExample5}
\Sel(E/{\bf Q}, E[3]) 
\simeq {\bf F}_{3} \oplus {\bf F}_{3} \oplus {\bf F}_{3}.
\end{equation}

We have a rational point $P$ of infinite order, so the rank of $E({\bf Q})$ is 
$\geq 1$. 
Take $N=3$ and consider $\ell=271$. 
Since $\tilde{\delta}_{271}=a_{1}^{(271)} = 35325 \equiv 9$ (mod $27$), 
$9$ is in   
$\Fitt_{1, {\bf Z}/p^{3}{\bf Z}}(\Sel(E/{\bf Q}, E[3^{3}])^{\vee})$ 
by Corollary \ref{HigherFittingIdeal2}. 
This implies that $\rank E({\bf Q})=1$ and 
$\#\sha(E/{\bf Q})[3^{\infty}] \leq 9$.
This together with (\ref{NExample5}) implies that 
\begin{equation} \label{NExample52}
\sha(E/{\bf Q})[3^{\infty}] \simeq {\bf Z}/3{\bf Z} \oplus {\bf Z}/3{\bf Z}.
\end{equation}

Note that if we used only Theorem \ref{MT2} and these computations, 
we could not get (\ref{NExample5}) nor (\ref{NExample52}) 
because we could not determine 
$\Theta_{1}({\bf Q})^{(\delta)}$ by finite numbers of computations. 
We need Theorem \ref{MMT2} to obtain (\ref{NExample5}) and (\ref{NExample52}).

\vspace{5mm}

\noindent {\bf (6)} 
For positive integers $d$ which are conductors of even Dirichlet 
characters (so $d=4m$ or $d=4m+1$ for some $m$)
satisfying $1 \leq d \leq 1000$, 
$d \equiv 1$ (mod $3$),
and $d \not \equiv 0$ (mod $11$), we computed 
$\Sel(E/{\bf Q}, E[3])$. 
Then $\dim \Sel(E/{\bf Q}, E[3])=0,1,2,3$, and 
the case of dimension $=3$ occurs only for $d=853$ in Example (5). 

\vspace{5mm}

\noindent {\bf (7)} We also considered negative twists. 
Take $d=-2963$. 
In this case, we know $L(E, 1) \neq 0$ and 
$L(E,1)/\Omega_{E}^{+}=81$.
We know from the main conjecture that 
the order of the $3$-component of 
$\sha(E/{\bf Q})$ is $81$, 
but {\it the main conjecture does not tell the structure of this group}.
Take $N=1$ and $\ell=19$. 
Then we compute $a_{2}^{(19)} = 2753/2 \not \equiv 0$ (mod $3$)
(we have $\vartheta_{{\bf Q}(19)} \equiv 
-432S+(2753/2)S^{2}$ mod $(9, S^{3})$). 
Since $c_{2}=1$, this shows that 
$a_{2}^{(19)}$ is in 
$\Fitt_{2, {\bf F}_{3}}(\Sel(E/{\bf Q}, E[3])^{\vee})$ 
by Theorem \ref{HigherFittingIdeal}. 
Therefore, we have 
$\Fitt_{2, {\bf F}_{3}}(\Sel(E/{\bf Q}, E[3]))={\bf F}_{3}$, which 
implies that 
$\Sel(E/{\bf Q}, E[3]) \simeq ({\bf F}_{3})^{\oplus 2}$. 
This denies the possibility of 
$\sha(E/{\bf Q})[3^{\infty}]  \simeq ({\bf Z}/3{\bf Z})^{\oplus 4}$, 
and we have  
$$\sha(E/{\bf Q})[3^{\infty}]  \simeq {\bf Z}/9{\bf Z} \oplus {\bf Z}/9{\bf Z}.$$

\vspace{5mm}

\noindent {\bf (8)} Let $E$ be the curve $y^{2}+xy+y=x^{3}+x^{2}-15x+16$ 
which is 563A1 in Cremona's book \cite{Cre}.
We take $p=3$. 
Since $a_{3} =-1$, $\Tam(E)=1$, $\mu=0$ and the Galois representation on 
$T_{3}(E)$ is surjective, all the conditions we assumed are satisfied. 
We know $\epsilon=1$ and $L(E,1)=0$.
Take $N=1$. 
We compute 
${\cal P}_{1}=\{13, 61, 103, 109, 127, 139,...\}$.
For admissible pairs $\{13,103\}$, $\{13,109\}$, 
we compute $\tilde{\delta}_{13 \times 103}=-6819 \equiv 0$ (mod $3$) 
and $\tilde{\delta}_{13 \times 109}=-242 \not \equiv 0$ (mod $3$). 
From the latter, we know that 
$$s_{13 \times 109}: 
\Sel(E/{\bf Q}, E[3]) \stackrel{\simeq}{\longrightarrow} 
({\bf F}_{3})^{\oplus 2}$$
is bijective by Theorem \ref{MMT2} (3).
Since $\lambda'=2$, the main conjecture also holds by 
Corollary \ref{10Prop3C}.
We know $L(E,1)=0$, so $\Sel(E/{\bf Q}, E[3^{\infty}]) 
\simeq ({\bf Z}_{3})^{\oplus 2}$. 

Numerically, we can find rational points 
$P=(2,-2)$ and $Q=(-4,7)$ on this elliptic curve. 
We can check that $E({\bf F}_{13})$ is cyclic of order $12$,  
$E({\bf F}_{103})$ is cyclic of order $84$, and 
$E({\bf F}_{109})$ is cyclic of order $102$.
The points $P$ and $Q$ have the same image and do not vanish in 
$E({\bf F}_{13})/3E({\bf F}_{13})$, 
but the image of $P$ in $E({\bf F}_{109})/3E({\bf F}_{109})$ 
is zero, and  the image of $Q$ in $E({\bf F}_{109})/3E({\bf F}_{109})$ 
is non-zero. 
This shows that $P$ and $Q$ are linearly independent over ${\bf Z}_{3}$, 
and $s_{13 \times 109}$ is certainly bijective.
Since all the elements in $\Sel(E/{\bf Q}, E[3^{\infty}])$ come from 
the points, we have 
$\sha(E/{\bf Q})[3^{\infty}]=0$.
On the other hand, 
the image of $P$ in $E({\bf F}_{103})/3E({\bf F}_{103})$ 
coincides with the image of $Q$, so 
$s_{13 \times 103}$ is not bijective. 
This is an example for which $\tilde{\delta}_{13 \times 103} \equiv 0$ 
(mod $3$) and $s_{13 \times 103}$ is not bijective.

\vspace{5mm}

\noindent {\bf (9)} Let $E$ be the elliptic curve 
$y^{2}+xy+y=x^{3}+x^{2}-10x+6$ 
which has conductor $18097$. 
We take $p=3$. 
We know $a_{3} =-1$, $\Tam(E)=1$, $\mu=0$ and the Galois representation on 
$T_{3}(E)$ is surjective, so all the conditions we assumed are satisfied. 
In this case, $\epsilon=-1$ and $L(E,1)=0$.
Take $N=1$. 
We compute 
${\cal P}_{1}=\{7, 19, 31, 43, 79,...,601,...\}$.
We know $\{7,43,601\}$ is admissible. 
We have $\tilde{\delta}_{7 \times 43 \times 601}=-2424748 
\not \equiv 0$ (mod $3$), and $7 \times 43 \times 601$ 
is $\delta$-minimal. 
We thank K. Matsuno heartily for his computing this value for us. 
The group $E({\bf F}_{7})$ is cyclic of order $9$ and $E({\bf F}_{43})$ 
is cyclic of order $42$. 
The point $(0,2)$ is on this elliptic curve, and has non-zero image both 
in $E({\bf F}_{7})/3E({\bf F}_{7})$ and 
$E({\bf F}_{43})/3E({\bf F}_{43})$. 
So both $s_{7}$ and $s_{43}$ are surjective, and 
we can apply Theorem \ref{MMT2} (4) to get 
$$s_{7 \times 43 \times 601}: 
\Sel(E/{\bf Q}, E[3]) \stackrel{\simeq}{\longrightarrow} 
({\bf F}_{3})^{\oplus 3}$$
is bijective. 

Numerically, we can find 3 rational points $P=(0,2)$, $Q=(2,-1)$, $R=(3,2)$ 
on this elliptic curve, and easily check that 
the restriction of $s_{7 \times 43 \times 601}$ to 
the subgroup generated by $P$, $Q$, $R$ in $\Sel(E/{\bf Q}, E[3])$ 
is surjective. 
Therefore, we have checked numerically that 
$s_{7 \times 43 \times 601}$ is bijective. 
This also implies that $\rank E({\bf Q})=3$ 
since $E({\bf Q})_{\tors}=0$.
Therefore, all the elements of $\Sel(E/{\bf Q}, E[3^{\infty}])$ come from 
the rational points, and we have $\sha(E/{\bf Q})[3^{\infty}]=0$. 

\subsection{A Remark on ideal class groups} \label{subsection54}

We consider the classical Stickelberger element 
$$\tilde{\theta}_{{\bf Q}(\mu_{m})}^{St}=
\sum_{\stackrel{\scriptstyle a=1}{(a,m)=1}}^{m} 
(\frac{1}{2}-\frac{a}{m}) \sigma_{a}^{-1} 
\in {\bf Q}[\Gal({\bf Q}(\mu_{m})/{\bf Q})]$$
(cf. (\ref{ModularElement})).
Let $K={\bf Q}(\sqrt{-d})$ be an imaginary quadratic field 
with conductor $d$, and $\chi$ be the corresponding quadratic character.
Let $m$ be a squarefree product whose prime divisors $\ell$ split in $K$ and 
satisfy $\ell \equiv 1$ (mod $p$).  
Using the above classical Stickelberger element, we define 
$\tilde{\delta}_{m,K}^{St}$ by 
$$\tilde{\delta}_{m,K}^{St}=-\sum_{\stackrel{\scriptstyle a=1}{(a,md)=1}}^{md} 
\frac{a}{md} \chi(a) 
(\prod_{\ell \mid m} \log_{{\bf F}_{\ell}}(a)) $$
(cf. (\ref{ModularSymbolDelta})).
We denote by $Cl_{K}$ the class group of $K$, and define 
the notion ``$\delta_{K}^{St}$-minimalness" analogously. 
We consider the analogue of Conjecture \ref{Conj2} for $\tilde{\delta}_{m,K}^{St}$ 
and $\dim_{{\bf F}_{p}}(Cl_{K}/p)$. 
Namely, we ask whether 
$\dim_{{\bf F}_{p}}(Cl_{K}/p)=\epsilon(m)$ for 
a $\delta_{K}^{St}$-minimal $m$.
Then the analogue does not hold. 
For example, take $K={\bf Q}(\sqrt{-23})$ and $p=3$. 
We know $Cl_{K} \simeq {\bf Z}/3{\bf Z}$. 
Put $\ell_{1}=151$ and $\ell_{2}=211$. 
We compute 
$\tilde{\delta}_{\ell_{1},K}^{St}=-270 \equiv 0$ (mod $3$), 
$\tilde{\delta}_{\ell_{2},K}^{St}=-1272 \equiv 0$ (mod $3$), 
and 
$\tilde{\delta}_{\ell_{1} \cdot \ell_{2},K}^{St}=
-415012 \equiv 2$ (mod $3$). 
This means that $\ell_{1} \cdot \ell_{2}$ is $\delta_{K}^{St}$-minimal. 
But, of course, we know 
$\dim_{{\bf F}_{p}}(Cl_{K}/p)=1 < 2=\epsilon(\ell_{1} \cdot \ell_{2})$.

{\small

\hspace{7cm} Masato {\sc Kurihara}

\hspace{7cm} Department of Mathematics, 

\hspace{7cm} Keio University, 

\hspace{7cm} 3-14-1 Hiyoshi, Kohoku-ku, 

\hspace{7cm} Yokohama, 223-8522, Japan

\hspace{7cm} {\tt kurihara@math.keio.ac.jp}

}


\begin{thebibliography}{99}
\bibitem{Cre} Cremona, J.E., {\it Algorithms for modular elliptic curves}, 
Cambridge University Press, 1992.
\bibitem{Gr1} Greenberg, R., Iwasawa theory for $p$-adic representations, 
in {\it Algebraic Number Theory -- in honor of K. Iwasawa}, 
Advanced Studies in Pure Math 17 (1989), 97-137.
\bibitem{Gr2} Greenberg, R., Iwasawa theory for elliptic curves, 
in {\it Arithmetic theory of elliptic curves}, Cetraro, Italy 1997, 
Springer Lecture Notes in Math 1716 (1999), 51-144.
\bibitem{Gr3} Greenberg, R., {\it Iwasawa theory, projective modules, 
and modular representation}, Memoirs of the AMS Number 992 (2011).
\bibitem{Grei3} Greither, C., Computing Fitting ideals of Iwasawa modules, 
Math. Zeitschrift 246 (2004), 733-767.
\bibitem{HM1} Hachimori, Y. and Matsuno K., An analogue of Kida's formula 
for the Selmer groups of elliptic curves, J. Algebraic Geometry 8 (1999), 581-601.
\bibitem{KK1} Kato, K., $p$-adic Hodge theory and values of zeta functions 
of modular forms, in {\it Cohomologies $p$-adiques et applications 
arithm\'{e}tiques III}, Ast\'{e}risque 295 (2004), 117-290.
\bibitem{Kol} Kolyvagin, V.A., Euler systems, The Grothendieck 
Festschrift Vol II (1990), 435-483.
\bibitem{Ku4} Kurihara, M., Iwasawa theory and Fitting ideals, 
J. reine angew. Math. 561 (2003), 39-86.
\bibitem{Ku3} Kurihara, M., On the structure of ideal class groups of 
CM-fields, Documenta Mathematica, Extra Volume Kato (2003), 
539-563. 
\bibitem{Ku5} Kurihara, M., Refined Iwasawa theory and 
Kolyvagin systems of Gauss sum type, Proceedings of the 
London Mathematical Society 104 (2012), 728-769. 
\bibitem{Ku6} Kurihara, M., Refined Iwasawa theory for $p$-adic representations 
and the structure of Selmer groups, to appear in M\"{u}nster Journal of Mathematics
http://www.math.keio.ac.jp/~kurihara/
\bibitem{HM2} Matsuno K., An analogue of Kida's formula for 
the $p$-adic $L$-functions of modular elliptic curves, J. Number Theory 84 
(2000), 80-92.
\bibitem{MR} Mazur, B. and Rubin, K, {\it Kolyvagin systems}, 
Memoirs of the AMS Vol 168, Number 799 (2004).
\bibitem{MR2} Mazur, B. and Rubin, K, Organizing the arithmetic of 
elliptic curves, Advances in Mathematics 198 (2005), 504-546.
\bibitem{MT} Mazur, B. and Tate J., Refined conjectures of the 
``Birch and Swinnerton-Dyer type", Duke Math. J. 54 (1987), 711-750.
\bibitem{Ne0} Nekov\'{a}\v{r}, J., On the parity of ranks 
of Selmer groups. II, C. R. Acad. Sci. Paris S\'{e}r. I Math. 332 (2001), 
99-104. 
\bibitem{Ne} Nekov\'{a}\v{r}, J., {\it Selmer complexes}, 
Ast\'{e}risque No. 310 (2006). 
\bibitem{N} Northcott, D. G., {\it Finite free resolutions}, 
Cambridge Univ. Press (1976).
\bibitem{Pope} Popescu, C. D., On the Coates-Sinnott conjecture, 
Mathematische Nachrichten 282 (2009), 1370-1390. 
\bibitem{Ru0} Rubin, K., The main conjecture, Appendix to 
{\it Cyclotomic fields I and II} by S. Lang, 
Graduate Texts in Math. 121, Springer-Verlag (1990), 397-419.
\bibitem{Ru1} Rubin, K., Kolyvagin's system of Gauss sums, 
{\it Arithmetic Algebraic Geometry}, G. van der Geer et al eds, 
Progress in Math 89 (1991) 309-324.
\bibitem{Ru2} Rubin, K., {\it Euler systems}, Annals of Math. Studies 147, 
Princeton Univ. Press 2000.
\bibitem{Sch0} Schneider, P., Iwasawa $L$-functions of varieties over 
algebraic number fields, A first approach, Invent. math. 71 (1983), 251-293.
\bibitem{CL} Serre, J.-P., {\it Corps Locaux}, Hermann, Paris 1968 
(troisi\`{e}me \'{e}dition).
\bibitem{SU1} Skinner, C. and Urban, E., The Iwasawa main conjecture for 
$GL_{2}$, Invent. math. 195 (2014), 1-277.
\bibitem{Stevens} Stevens, G., Stickelberger elements and modular parametrizations 
of elliptic curves, Invent math 98 (1989), 75-106.
\end{thebibliography}
\end{document}